\documentclass[a4paper,10pt,english]{smfart}
\usepackage{amsfonts}
\usepackage{amsmath} 
\usepackage{hyperref} 

\numberwithin{equation}{section}
\usepackage{latexsym}
\usepackage{array}
\usepackage{amssymb}
\usepackage{smfthm}
\usepackage{bull}
\theoremstyle{plain}

\newtheorem*{acknowledgements}{Acknowledgements}
\newcommand{\R}{  \mathbb{R}   }
\newcommand{\eps}{\varepsilon}
\newcommand{\e}{  \text{e}   }
\newcommand{\C}{  \mathbb{C}   }

\newcommand{\N}{  \mathbb{N}   }

\newcommand{\h}{  \hbar   }

\author{Laurent Thomann}
\address{Universit\'e Paris-Sud, Math\'ematiques, B\^at 425\\Tel 0169155785\\ 91405
Orsay Cedex.}
\email{ laurent.thomann@math.u-psud.fr}
\urladdr{http://www.math.u-psud.fr/~thomann}
\title[Instabilities for supercritical NLS]{ Instabilities for
  supercritical Schr\"odinger equations in analytic manifolds}

\begin{document}

\frontmatter
 \begin{abstract}
In this paper we consider supercritical nonlinear Schr\"odinger
equations in an analytic Riemannian manifold $(M^d,g)$, where the
metric $g$ is analytic. Using an analytic WKB method, we are able to
construct an Ansatz for the semiclassical equation for times independent of the
small parameter. These approximate solutions will help to show two
different types of instabilities. The first is in the energy space, and
the second is an immediate loss of regularity in higher Sobolev norms. 
\end{abstract}

\maketitle
\mainmatter


\section{Introduction}

 \indent Let $(M^d,g)$ be an analytic Riemannian manifold of dimension
 $d\geq3$. In all the paper we assume that the metric $g$ is
 analytic. Let $p$ an odd integer. \\
 We consider the nonlinear Schr\"odinger equation

\begin{equation}\label{nls}
\left\{
\begin{aligned}
&i \partial_t u+\Delta_{g} u   =  \omega|u|^{p-1} u,\quad
(t,x)\in\R\times M^d,\\
&u(0,x)= u_0(x),
\end{aligned}
\right.
\end{equation}
with  either $\omega=1$ (defocusing equation) or $\omega=-1$ (focusing
equation).\\
Here $\Delta=\Delta_g$ denotes the Laplace-Beltrami operator defined
by $\Delta=\text{div} \nabla$.\\
It is known that the mass
\begin{equation}\label{L2}
\|u(t)\|_{L^2( M^d)}=\|u_0\|_{L^2( M^d)},
\end{equation}
and the energy
\begin{equation}\label{Hamil}
H(u)(t)=\int_{ M^d}\Big(\frac12|\nabla u|^2+\frac{\omega}{p+1}|u|^{p+1}\Big)\text{d}x=H(u_0,\omega),
\end{equation}
are conserved by the flow of \eqref{nls}, at least formally.\\
Denote also by
\begin{equation}\label{Hamil+}
H^+(u)=\int_{ M^d}\Big(\frac12|\nabla u|^2+\frac1{p+1}|u|^{p+1}\Big)\text{d}x.
\end{equation}

\noindent In the following we will need the definition of uniform
well-posedness :
\begin{defi}
Let $X$ be a Banach space. We say that the Cauchy problem \eqref{nls}
is locally uniformly well-posed in $X$, if for any bounded subset
$\mathcal{B}\subset X$, there exists $T>0$ and a solution
$u\in\mathcal{C}\big([-T,T];X\big)$ of \eqref{nls} and such that the
  flow map 
$$u_0 \in \mathcal{B}  \longmapsto u(t)=\Phi_t(u_0)\in X,$$
is uniformly continuous for any $-T\leq t\leq T$.
\end{defi}
\subsection{Instability in the energy space}
$~$\\[5pt]
By the works of J. Ginibre and G. Velo \cite{GV}, T. Cazenave and
F. B.  Weissler \cite{CW},  we know that \eqref{nls} is
locally uniformly well-posed in the energy space $X= H^1(\R^d)\cap
L^{p+1}(\R^d) = H^1(\R^d)$ when $p< (d+2)/(d-2)$.\\
Our first result states that this result does not hold when $p>
(d+2)/(d-2)$ is an odd integer.

\begin{theo}\label{thm1}
Let $p>(d+2)/(d-2)$ be an odd integer, $\omega \in \{-1,1\}$, and let $H^+$ be given by
\eqref{Hamil+}. Let $m\in M^d$. There exist a positive  sequence $r_n\longrightarrow
0$, and two  sequences $u_0^n, \widetilde{u_0}^n\in
\mathcal{C}_0^{\infty}( M^d)$ of  Cauchy data
    with support in the ball $\big\{|x-m|_g\leq r_n\big\}$,  a sequence of times
$t_n\longrightarrow 0$, and  constants $c,C>0$ such that
\begin{equation}\label{st1}
H^+(u_0^n)\leq C, \quad H^+(\widetilde{u_0}^n)\leq C,
\end{equation}
\begin{equation}\label{st2}
H^+(u_0^n-\widetilde{u_0}^n)\longrightarrow 0,\quad \text{when} \quad
n\longrightarrow +\infty,
\end{equation}
and such that the solutions $u^n$, $\tilde{u}^n$  of \eqref{nls} satisfy
\begin{equation}\label{st3}
\limsup_{n\to +\infty}
\int_{ M^d}\big|(u^n-\widetilde{u}^n)(t_n)\big|^{p+1}\text{d}x>c.
\end{equation}
Moreover, the sequences $u_0^n, \widetilde{u_0}^n$ can be chosen such
that there exist $\nu_0>0$ and $q_0>p+1$, such that for all
$0\leq \nu<\nu_0$ and $p+1\leq q<q_0$,
\begin{equation}\label{st4}
\|u_0^n-\widetilde{u_0}^n\|_{H^{1+\nu}( M^d)}+
\|u_0^n-\widetilde{u_0}^n\|_{L^{q}( M^d)}\longrightarrow 0.
\end{equation}

\end{theo}
\noindent For $k\in \R$, the norm $\|\cdot\|_{H^k(M^d)}$ is defined by 
$$ \|f\|_{H^k(M^d)} =\|(1-\Delta)^{k/2}f\|_{L^2(M^d)}.$$ 
R. Carles \cite{Carles2} obtains a similar result for the defocusing
cubic equation in $\R^d$. 
\noindent An analog of Theorem \ref{thm1} was proved by G. Lebeau
\cite{Lebeau} for the supercritical wave equation, but for a
nonlinearity of the form $u^p$. After a rescaling of \eqref{nls} to a semiclassical equation,
we also have an almost  finite speed of propagation principle. This is one 
reason why such a result was expected for nonlinear supercritical
Schr\"odinger equations.
\subsection{Ill-posedness in Sobolev spaces}
$~$\\[5pt]
\noindent Assume here that $(M^d,g)$ is the euclidian space with the
    canonical metric $(M^d,g) =(\R^d,\text{can})$. Let $T>0$ and let $u:]-T,T[\times \R^d\longrightarrow \C$
    satisfy \eqref{nls}. Then for all $\lambda\in\R$
\begin{eqnarray}
\begin{array}{ccc}
u^{\lambda}\colon]-\lambda^{-2}T,\lambda^{-2}T[\times\R^d &\longrightarrow & \C \nonumber\\[5pt]
(t,x)&\longmapsto & u^{\lambda}(t,x)=\lambda^{\frac{2}{p-1}
    }u(\lambda^2t,\lambda x),
\end{array}
\end{eqnarray}
is also a solution of \eqref{nls}.\\
Define the critical index for Sobolev well-posedness 
\begin{equation}\label{defsigmacrit}
\sigma_{\text{c}}=\frac{d}2-\frac{2}{p-1}.
\end{equation}
Then, for all $f\in\dot{H}^{\sigma_{\text{c}}}(\R^d)$ (the homogeneous Sobolev space)
and $\lambda \in \R$
$$\lambda^{\frac{2}{p-1}}\|f(\lambda\cdot)\|_{\dot{H}^{\sigma_{\text{c}}}(\R^d)}=\|f\|_{\dot{H}^{\sigma_{\text{c}}}(\R^d)}.$$
This scaling notion is relevant, as we have the following results :\\[5pt]
$\bullet$ Let $\sigma> \sigma_{\text{c}}$, then the equation \eqref{nls} is
locally uniformly well-posed in $X={H}^{\sigma}(\R^d)$, \cite{GV},\cite{CW}.\\[5pt]
$\bullet$ If $0<\sigma<\sigma_{\text{c}}$, the problem \eqref{nls} is
ill-posed in ${H}^{\sigma}(\R^d)$, in the sense that there exist a sequence of initial data
$u^n_0$ so that 
$$\|u^n_0\|_{H^{\sigma}(\R^d)}\longrightarrow 0,$$
and a sequence of times $t_n\longrightarrow 0$ such that the solution
$u^n$ of \eqref{nls} satisfies 
$$\|u^n(t_n)\|_{H^{\rho}(\R^d)}\longrightarrow +\infty,$$
for $\rho=\sigma$ (see Christ-Colliander-Tao \cite{CCT2}), or even for all
$\rho\in]\sigma/(\frac{d}2-\sigma),\sigma]$ in the particular case $\omega=1$ and
    $p=3$, (see Carles \cite{Carles} and Alazard-Carles \cite{AlCa}).\\[10pt]
Here we prove

\begin{theo}\label{thm2}
Assume that $(M^d,g) =(\R^d,\text{can})$. Let $p\geq 3$ be an odd integer, $\omega \in \{-1,1\}$, and let
$0<\sigma<d/2-2/{(p-1)}$. There
exist a  sequence  $\check{u}_0^n\in
\mathcal{C}^{\infty}(\R^d)$ of  Cauchy data and a sequence of times
$\tau_n\longrightarrow 0$ such that
\begin{equation}\label{1.6}
\|\check{u}_0^n\|_{H^{\sigma}(\R^d)}\longrightarrow 0,\quad \text{when} \quad
n\longrightarrow +\infty,
\end{equation}
and such that the solution $\check{u}^n$ of \eqref{nls} satisfies
\begin{equation}\label{1.7}
\|\check{u}^n(\tau_n)\|_{H^{\rho}(\R^d)}\longrightarrow +\infty,\quad \text{when} \quad
n\longrightarrow +\infty,\quad \text{for all} \quad \rho\in\Big]\frac{\sigma}{\frac{p-1}{2}(\frac{d}2-\sigma)},\sigma\Big].
\end{equation}
\end{theo}

\noindent In the general case of an analytic manifold  $(M^d,g)$ with an
analytic metric $g$, we obtain the weaker result 

\begin{theo}\label{thm3}
Let $p\geq 3$ be an odd integer, $\omega \in \{-1,1\}$, and let
$0<\sigma<d/2-2/{(p-1)}$. Let $m\in M^d$. There
exist a positive sequence $r_n\longrightarrow 0$ and a sequence  $\check{u}_0^n\in
\mathcal{C}_0^{\infty}(M^d)$ of  Cauchy data
    with support in the ball $\big\{|x-m|_g\leq r_n\big\}$,  a sequence of times
$\tau_n\longrightarrow 0$ such that
\begin{equation*}
\|\check{u}_0^n\|_{H^{\sigma}(M^d)}\longrightarrow 0,\quad \text{when} \quad
n\longrightarrow +\infty,
\end{equation*}
and such that the solution $\check{u}^n$ of \eqref{nls} satisfies
\begin{equation*}
\|\check{u}^n(\tau_n)\|_{H^{\rho}(M^d)}\longrightarrow +\infty,\quad \text{when} \quad
n\longrightarrow +\infty,\quad \text{for all} \quad
\rho\in\Big]I(\sigma),\sigma\Big],
\end{equation*}
where $I(\sigma)$ is defined by
\begin{equation*}
I(\sigma)=\left\{\begin{array}{ll} 
\frac{\sigma}{2} \quad &\text{for} \quad 0<\sigma\leq\frac{d}{2}-\frac4{p-1}, \\[5pt]  
 \frac{\sigma}{\frac{p-1}{2}(\frac{d}2-\sigma)} \quad \quad &\text{for} \quad
\frac{d}{2}-\frac{4}{p-1}\leq\sigma<\frac{d}{2}-\frac{2}{p-1}.
\end{array} \right.
\end{equation*}
\end{theo}

\noindent In the case  $p=3$ and $\omega=1$, Theorem \ref{thm2} was shown by
R. Carles \cite{Carles} using the convergence of the WKB method for
$\mathcal{C^{\infty}}$ data (see \cite{PGX}, \cite{Grenier} ). Recently, T. Alazard and R. Carles \cite{AlCa2} have obtained a justification for nonlinear geometric optics when $p>3$ with $H^{\infty}$ data.\\
Consider now the semiclassical equation
\begin{equation}\label{eq0}
ih \partial_t v+h^2\Delta v   =  |v|^{p-1} v.
\end{equation}
In \cite{AlCa}, T. Alazard and R. Carles prove that for all non
trivial initial condition $v(0,\cdot)\in \mathcal{S}(\R^d)$, the solution $v$ of
\eqref{eq0} oscillates immediately: There exists  $\tau>0$ so that 
$$ \liminf_{h\to0}\||h\nabla|^s v(\tau) \|_{L^2(\R^d)}>0,$$
for all $s\in]0,1]$.
This yields the result of Theorem \ref{thm2} for the defocusing
equation in the euclidien space for any smooth Cauchy
condition. Their method does not apply to the focusing case.\\[10pt] 
Denote by $\sigma_{\text{sob}}$ the Sobolev exponent so that
$ {\dot{H}}^{\sigma_{\text{sob}}}(\R^d)\subset L^{p+1}(\R^d)$, i.e.
\begin{equation}\label{defsigmasob}
\sigma_{\text{sob}}=\frac{d}2-\frac{d}{p+1}.
\end{equation}
Let $p>(d+2)/(d-2)$, then $\sigma_{\text{sob}}<\sigma_{\text{c}}$. As pointed out by
G. Lebeau and R. Carles, for $\sigma=\sigma_{\text{sob}}$, Theorem \ref{thm2} yields 
$$\|\check{u}^n_0\|_{H^{\sigma_{\text{sob}}}(\R^d)}\longrightarrow 0,\quad 
\|\check{u}^n(\tau_n)\|_{H^{\rho}(\R^d)}\longrightarrow +\infty, $$
for all $\rho\in]1,\sigma_{\text{sob}}]$. This interval can not be
    enlarged. Indeed, for all $\rho\leq 1$, the conservation of the quantities \eqref{L2} and
    \eqref{Hamil} together with the embedding 
   $ \dot{H}^{\sigma_{\text{sob}}}(\R^d)\subset L^{p+1}(\R^d)$ yield for all $\tau>0$
$$\|\check{u}^n(\tau)\|_{H^{\rho}(\R^d)}\longrightarrow 0. $$
See also \cite{BGI}.\\[5pt]
G. Lebeau \cite{Lebeau2} obtains a stronger result for the wave
equation in $(\R^d,\text{can})$, with the
same range for $\rho$ in \eqref{1.7}, but the loss of derivatives is
obtained with only one one Cauchy condition, instead of a sequence. \\[5pt]
Theorem \ref{thm1} can not be deduced from Theorem \ref{thm2}. In
fact, the sequences constructed with $\sigma=1$ such that 
$$\|\check{u}^n_0\|_{H^{1}(M^d)}\longrightarrow 0,\quad 
\|\check{u}^n(\tau_n)\|_{H^{1}(M^d)}\longrightarrow +\infty, $$
satisfy $H^+(\check{u}^n_0)\longrightarrow +\infty$, when $n$ tends to
infinity.\\[10pt]
The instabilities of Theorems \ref{thm1},  \ref{thm2} and  \ref{thm3} are not
geometrical effects, they are only caused by the high exponent of the
nonlinearity.\\
We could also consider more general analytic nonlinearities, for
instance\\ $\pm (1+|u|^2)^{\alpha/2}u$ with $\alpha>(d+2)/(d-2)$.\\
Notice that the focusing case with non analytic Cauchy conditions is
more intricate, as other phenomenons are involved, like finite time
explosion.\\[5pt] 
The main ingredient of the proof of our results is the construction of
approximate solutions of \eqref{nls}, via analytic nonlinear geometric
optics, as done by P. G\'erard in \cite{PGX}. This work will be
adapted to the case $(M^d,g)=(\R^d,\text{can})$. We will work in
weighted spaces, so that these solutions concentrate in a point of
$\R^d$, and then the construction in $(M^d,g)$ will follow directly, as
we are able to work only in one local chart.
\\[10pt]
The plan of the paper is the following\\[5pt]
1. We first construct a formal solution of \eqref{nls}.\\
 \indent  a) In Section \ref{SectionWKB} we deal with the case
 $(M^d,g)=(\R^d,\text{can})$ : First we reduce \eqref{nls} to a semiclassical equation as done in
\cite{Lebeau} and \cite{Carles}, then we adapt the analytic WKB method given
in \cite{PGX} to $\R^d$.\\
\indent b) In Section \ref{SectionAna} we consider the general case
of an analytic manifold with an analytic metric.\\
2. We obtain a family of approximate solutions of
\eqref{nls}. (Section \ref{SectionAnsatz})\\
3. Using two different rates of   concentration of this family, we
prove the main results. (Section \ref{SectionInstab})

\begin{enonce}{Notations}
In this paper $c$, $C$ denote constants the value of which may change
from line to line. These constants will always be independent of $h$. We use the notations $a\sim b$,  
$a\lesssim b$ if $\frac1C b\leq a\leq Cb$ , $a\leq Cb$
respectively. We write $a\ll b$ if $a\leq Kb$ for some large constant
$K$ which is  independent of $h$.
\end{enonce}

\begin{acknowledgements}
The author would like to thank N. Burq his adviser for this
interesting subject and his guidance, and P. G\'erard for giving his permission to
reproduce a part of the work \cite{PGX} in the appendix. The author is
also grateful to S. Alinhac and T. Alazard for many enriching
discussions and clarifications.
\end{acknowledgements}
\section{Nonlinear geometric optics }\label{SectionWKB}

\subsection{The Euclidian case}
\subsubsection{Reduction to a semiclassical equation}

\noindent Following \cite{Lebeau}, \cite{Carles}, we reduce the equation \eqref{nls} to a semiclassical
equation, and therefore make the following change of variables and
unknown function

\begin{equation}\label{chgt}
\left\{
\begin{aligned}
&t=\h^{\alpha}s,\quad x=\h z,\quad h=\h^{\beta},\\
&u(\h^{\alpha}s,\h z)=\h^{\gamma}v(s,z,h),
\end{aligned}
\right.
\end{equation}
where $h\in ]0,1]$ is a small parameter, and where $\beta>0$. The
    value of $\beta$ will be given in Section \ref{SectionInstab}, in
    terms of $p$ and $d$ to prove Theorem \ref{thm1}, and in terms of
    $p$, $d$ and $\sigma$ to  prove Theorem \ref{thm2}.\\
    If we choose 
\begin{equation}\label{parameters}
\alpha=\beta+2, \quad (p-1)\gamma=-2(\beta+1),
\end{equation}
we are lead to studying the Cauchy problem
\begin{equation}\label{nlsh}
\left\{
\begin{aligned}
&i h\partial_s v(s,z)+h^2\Delta v(s,z)   = \omega |v|^{p-1} v(s,z),\\
&v(0,z)= v_0(z).
\end{aligned}
\right.
\end{equation}

\noindent Following the ideas of nonlinear geometric optics, we can search a
solution of \eqref{nlsh} for small times (but independent of $h$) of the form
\begin{equation}\label{ansatz}
v(s,z,h)= a(s,z,h)\e^{iS(s,z)/h},
\end{equation}
where formally

\begin{equation}\label{ansatz2}
a(s,z,h) =\sum_{j\geq 0} a_j(s,z)h^j.
\end{equation}
Then $v$ is a formal solution of equation \eqref{nlsh} if the couple
$(S,a)$ satisfies  the system

\begin{equation}\label{systgeom}
\left\{
\begin{aligned}
&\partial_s S+(\nabla S)^2+\omega|a_0|^{p-1}=0,\\
&\partial_s a+2\nabla S\cdot\nabla
a+a\Delta S-ih\Delta a +\frac{i\omega a}{h}(|a|^{p-1}-|a_0|^{p-1})=0,\\
&S(0,z)=S^0(z),\;a(0,z,h)=a^0(z,h),
\end{aligned}
\right.
\end{equation}
where $v(0,z,h)=a^0(z,h)\e^{iS^0(z)/h}$.\\
In fact to obtain the system \eqref{systgeom}, plug \eqref{ansatz} in equation \eqref{nlsh} and identify the
coefficients in the expansion in 
powers of $h$. The first equation of \eqref{systgeom} corresponds to  the
coefficients of $h^0$, and the second to the others, after division by
$h$.
Notice that $S$ will be a real function, if the data
$S(0,\cdot)$ is real.\\
The WKB method consits now in plugging the developement given by
\eqref{ansatz2} in \eqref{systgeom}. Annihilating the coefficients of $h^j$, for $j\geq 0$,
yields a cascade of equations. And if we are able to solve them, this
gives an approximate solution $v_{\text{app}}$ of \eqref{nlsh}
\begin{equation}\label{WKB}
ih \partial_s v_{\text{app}}+h^2\Delta v_{\text{app}}   =  |v_{\text{app}}|^{p-1} v_{\text{app}}+\mathcal{O}(h^{\infty}).
\end{equation}
Unfortunately, the obtained system is not closed: the equation
which gives  $a_j$ depends on $a_{j+1}$.\\
Moreover, in general,  using \eqref{WKB}, we can show that  $v_{\text{app}}$ is
close to a solution of \eqref{nlsh} only for times
$s\in[0,Ch\log{\frac1{h}}]$. See \cite{PGX}, Corollaire 1.\\[10pt]
To obtain an Ansatz for $h$-independent times, we work in an analytic
frame. Thus in the following we will consider $z$ as a complex variable.


\noindent \subsubsection{Construction of a formal solution of \eqref{nlsh}}\label{SectAnsatz}

Here we adapt step by step the proof of P. G\'erard \cite{PGX} given
in the case of the torus $\mathbb{T}^d$ to the case $\R^d$.\\
We need Sj\"ostrand's definition \cite{Sj} of an analytic symbol.

\begin{defi}\label{def1}
We say that the formal series $b(s,z,h)=\sum_{j\geq 0}b_j(s,z)h^j$ is an analytic symbol
if there exist positive constants $s_0,l,A,B>0$ such that for all $j\geq 0$
\begin{equation}
(s,z)\mapsto b_j(s,z) \;\text{is an holomorphic function on}\;\; \{|s|<s_0\}\times
 \{|\text{Im}\;z|<l\} ,
\end{equation}
and
\begin{equation}\label{analb}
| b_j(s,z)|\leq AB^jj!\;\;\text{on}\;\;\{|s|<s_0\}\times
 \{|\text{Im}\;z|<l\}.
\end{equation}
\end{defi}
\noindent Notice that $b$ has to be analytic in both variables, $s$ and $z$.\\[10pt]
To obtain proper estimates in Sobolev norms later, we want to make
sure  that the
functions are small at infinity in the space variable. Therefore we
define the weight 
\begin{equation}\label{weight}
W(z)=\e^{(1+z^2)^{1/2}},
\end{equation}
where $z^2=z_1^2+\cdots +z_d^2 $ for any $z=(z_1,\cdots,z_d)\in\C^d$.
Notice that $W$ is analytic in the band  $\{|\text{Im}\;z|<\frac12\}$,
thus in the following we fix $l<\frac12$.\\[5pt]
We introduce the space $\mathcal{H}(s_0,l,B)$ composed of the analytic
symbols\\ $b=\sum_{j\geq 0}b_jh^j$ satisfying: there exist $A,B>0$ so
that  
\begin{equation}\label{expanalb}
|W(z)b_j(s,z)|\leq AB^jj!\;\;\text{on}\;\;\{|s|<s_0\}\times
 \{|\text{Im}\;z|<l\},\;\forall j\geq 0.
\end{equation}

\begin{equation}\label{defH}
\mathcal{H}(s_0,l,B)=\left\{\begin{array}{c} 
b=\sum_{j\geq 0}b_jh^j\;\text{is an analytic symbol on} \\[5pt]   
\big(\{|s|<s_0\}\times
 \{|\text{Im}\;z|<l\}\big)\;\text{s.t.}\; b_j \;\text{satisfies \eqref{expanalb}} 
\end{array} \right\}.
\end{equation}
Let $\eps<1/B$. For $0\leq \theta \leq 1$, we can endow $\mathcal{H}(s_0,l,B)$ with the norms 
\begin{equation*}
\|b\|_{\theta}=\sum_{j\geq 0}\frac{\eps^j}{j!}\sup_{0<\tau<1}\;
\sup_{|s|<s_0(1-\tau)}\;\sup_{|\text{Im}\;z|<l\tau}|W(z)b_j(s,z)|\Big(1-\tau-\frac{|s|}{s_0}\Big)^{j+\theta}.
\end{equation*}

\noindent Each of these norms makes   $\mathcal{H}(s_0,l,B)$ a complete
space.\\
In the following, fix $0<\eps<1/B$, and let $0<h<\eps$. Fix also
$s_0$, $B>0$ and $l<\frac{1}2$. Denote by
\begin{equation*}
\mathcal{H}=\mathcal{H}(s_0,l,B),
\end{equation*}
and define
 \begin{equation*}
\mathcal{H}^0=\mathcal{H}(0,l,B),
\end{equation*}
 the restriction to
$s=0$ of  $\mathcal{H}$, endowed with the induced norms. This is the
space of the initial conditions.\\
We will solve the system
\eqref{systgeom} in $(\mathcal{H},\|\cdot\|_1)$ with a fixed point
argument. The choice of the space and norms are inspired by abstract
versions of the Cauchy-Kowaleski theorem \cite{BaGou}.\\[10pt]
We first give some properties of these norms.
\begin{lemm}\label{lem1}
There exists $C>0$ such that for all $\theta  \in [0,1]$ and  $b^1,b^2\in\mathcal{H}$
\begin{equation}\label{estnorm}
\|b^1\,b^2\|_{\theta}\leq C\|b^1\|_0\|b^2\|_{\theta}.
\end{equation}
\end{lemm}

\begin{proof} Set 
$$\Omega=\big{\{}(\tau,s,z)\;|
\;0<\tau<1,\;|s|<s_0(1-\tau),\;|\text{Im}\;z|<l\tau \big{\}},$$ 
and denote by
\begin{equation*}
\sup_{\Omega}=\sup_{0<\tau<1}\;
\sup_{|s|<s_0(1-\tau)}\;\sup_{|\text{Im}\;z|<l\tau}.
\end{equation*}
Let 
\begin{equation*}
b^1=\sum_{j\geq0}b^1_j\,h^j,\quad\text{and}\quad
b^2=\sum_{j\geq0}b^2_j\,h^j,
\end{equation*}
be two elements of $\mathcal{H}$, then $b^1\,b^2$ can be written
\begin{equation}\label{produit}
b^1\,b^2=\sum_{j=0}^{\infty}\big(\sum_{k=0}^jb_k^{1}\,b_{j-k}^{2}\big)h^j.
\end{equation}
It is easy to check that  there exists $C>0$ so that 
\begin{equation}\label{exp}
|W(z)|\leq C|W(z) |^2,
\end{equation}
on ${|\text{Im}\;z|<\frac12}$.
Therefore by \eqref{produit} and \eqref{exp}
\begin{eqnarray*}
\|b^1\,b^2\|_{\theta}&=&\sum_{j=0}^{\infty}\frac{\eps^j}{j!}\sup_{\Omega}\big|W(z)\sum_{k=0}^jb_k^{1}\,b_{j-k}^{2}(s,z)\big|\big(1-\tau-\frac{|s|}{s_0}\big)^{j+\theta}\\
&\leq&C\sum_{k=0}^{\infty}\sum_{j=k}^{\infty}\frac{\eps^k}{k!}\sup_{\Omega}|W(z)b^1_k(s,z)|\big(1-\tau-\frac{|s|}{s_0}\big)^k\\
&&\quad\quad\quad\cdot\frac{\eps^{j-k}}{(j-k)!}\sup_{\Omega}|W(z)b^2_{j-k}(s,z)|\big(1-\tau-\frac{|s|}{s_0}\big)^{j-k+\theta}\\
&=&\|b^1\|_0\|b^2\|_{\theta}.
\end{eqnarray*}
\end{proof}
For $|s|<s_0$, denote by $\partial_s^{-1}$ the operator defined by 
\begin{equation}\label{operateur}
\partial_s^{-1}b=\int_0^{s}b(\sigma)\text{d}\sigma\quad\text{for}\quad b\in\mathcal{H},
\end{equation}
and  $\partial_s^{-2}=\partial_s^{-1}\circ \partial_s^{-1}$.
We then have the following 
\begin{lemm}\label{lem2}
i) Let $A$ be one of the operators 
\begin{equation*} 
b\mapsto\nabla_z b,\;b\mapsto h\Delta_z b,\;b\mapsto\frac1h(b-b_0),
\end{equation*}
then there exists $C>0$ such that for all $h\in]0,\eps[$ and $b\in \mathcal{H}$
\begin{equation}\label{inegalite1} 
\|\partial_s^{-1}Ab\|_1\leq Cs_0\|b\|_1.
\end{equation}
ii) For all $\theta \in ]0,1]$, there  exists $C_{\theta}$ such that
    for all $h\in]0,\eps[$ and $b\in \mathcal{H}$
\begin{equation} \label{inegalite2} 
\|\partial_s^{-1}b\|_{\theta}\leq C_{\theta}s_0\|b\|_1.
\end{equation}
iii) There exists $C>0$ such that for all $h\in]0,\eps[$ and  $b\in \mathcal{H}$
\begin{equation} \label{inegalite3} 
\|\partial_s^{-2}b\|_0\leq Cs_0\|b\|_1.
\end{equation}
\end{lemm}

\begin{proof}  We can assume that $\|b\|_1=1$. Then there exists a
nonnegative sequence $d=(d_j)_{j\geq0}$ satisfying $\sum_{j\geq
  0}d_j=1$ so that, for all $j\geq 0$, for all $0<\tau <1$,
  $|s|<s_0(1-\tau)$, $|\text{Im}\;z|<l\tau$, we have
\begin{equation}\label{ineg1}
|W(z)\,b_j(s,z)|\leq C\frac{j!}{\eps^j}\frac{d_j}{(1-\tau-\frac{|s|}{s_0})^{j+1}}.
\end{equation}
Denote by $\nabla=\nabla_z$.\\
Proof of $i)$\\
$\bullet$ We prove the inequality $\|\partial_s^{-1}\nabla b\|_1\leq
Cs_0\|b\|_1$. Let $0<\tau<1$, $|s|<s_0(1-\tau)$ and $|\text{Im}\;z|<l\tau$. Let  $\tau<\tau'<1$. By the Cauchy formula we deduce
that for all $|s'|\leq|s|$ and  $|\text{Im}\;z|<l\tau$
\begin{equation*}
|\nabla b_j(s',z)|\leq \frac{C}{\tau'-\tau}\sup_{|\text{Im}\;z'|<l\tau'}| b_j(s',z')|.
\end{equation*}
Thus, as $|\nabla W|\leq |W|$, for all $|s'|\leq|s|$ and $|\text{Im}\;z|<l\tau$
\begin{equation}\label{ineg2}
|W(z)\nabla b_j(s',z)|\leq \frac{C}{\tau'-\tau}\sup_{|\text{Im}\;z'|<l\tau'}|W(z')\, b_j(s',z')|.
\end{equation}
Then by \eqref{ineg1} and \eqref{ineg2} we obtain

\begin{equation*}
\Big|W(z)\int_0^{s}\nabla b_j(s',z)\text{d}s'\Big|\leq
C\frac{j!}{\eps^j}d_j
\int_0^{|s|}\frac{1}{\tau'-\tau}\frac{\text{d}|s'|}{(1-\tau'-\frac{|s'|}{s_0})^{j+1}}.
\end{equation*}
We now make the choice 
\begin{equation}\label{tauprime}
\tau'-\tau=1-\tau'-\frac{|s'|}{s_0},\quad\text{i.e.}\quad \tau'=\frac12(1+\tau-\frac{|s'|}{s_0}),
\end{equation}
then $\tau'$ satisfies $\tau<\tau'<1$ because $0<\tau<1$ and
$|s'|<(1-\tau)s_0$.\\ Moreover, \eqref{tauprime} yields
$$1-\tau'-\frac{|s'|}{s_0}=\frac12(1-\tau-\frac{|s'|}{s_0}), $$
therefore 
\begin{eqnarray}
\Big|W(z)\int_0^{s}\nabla b_j(s',z)\text{d}s'\Big|&\leq&
C\frac{j!}{\eps^j}d_j\int_0^{|s|}\frac{\text{d}|s'|}{(1-\tau-\frac{|s'|}{s_0})^{j+2}}\nonumber\\
&\leq&Cs_0\frac{j!}{\eps^j}d_j\big((1-\tau-\frac{|s|}{s_0})^{-j-1}-(1-\tau)^{-j-1}\big).\nonumber
\end{eqnarray}
And thus,  as $|s|<s_0(1-\tau)$
\begin{eqnarray}
\frac{\eps^j}{j!}\Big|W(z)\int_0^{s}\nabla
b_j(s',z)\text{d}s'\Big|\Big(1-\tau-\frac{|s|}{s_0}\Big)^{j+1}&\leq&
Cs_0\Big((1-(1-\frac{|s|}{s_0(1-\tau)})^{j+1}\Big)d_j\nonumber\\
&\leq& Cs_0d_j.\nonumber
\end{eqnarray}
Finally, by the previous inequality
\begin{eqnarray}
\|\partial_s^{-1}\nabla b\|_1&=&\sum_{j\geq 0}\frac{\eps^j}{j!}\sup_{0<\tau<1}\;
\sup_{|s|<s_0(1-\tau)}\;\sup_{|\text{Im}\;z|<l\tau}\big|W(z)\int_0^{s}\nabla
b_j(s',z)\text{d}s'\big|\big(1-\tau-\frac{|s|}{s_0}\big)^{j+1}\nonumber\\
&\leq&
Cs_0\sum_{j\geq 0}d_j\leq Cs_0,\nonumber
\end{eqnarray}
which was the claim.\\[5pt]
$\bullet$
The inequality $h\|\partial_s^{-1}\Delta b\|_1\leq
Cs_0\|b\|_1$ can be shown by the same manner, using that $h<\eps$
compensates the loss of one more derivative.\\[5pt]
$\bullet$ 
Denote by $b'=(b-b_0)/h$, then for all $j\geq 0$, $b'_j=b_{j+1}$. By \eqref{ineg1}
\begin{equation*}
\Big|W(z)\int_0^{s} b_{j+1}(s',z)\text{d}s'\Big|\leq \frac{s_0}{j+1}\frac{(j+1)!}{\eps^{j+1}}d_{j+1}\big((1-\tau-\frac{|s|}{s_0})^{-j-1}-(1-\tau)^{-j-1}\big),
\end{equation*}
and therefore, for all $0<\tau <1$, $|s|<s_0(1-\tau)$, $|\text{Im}\;z|<l\tau$ and
$j\geq0$
\begin{eqnarray}
\frac{\eps^j}{j!}\Big|W(z)\int_0^{s}
b_{j+1}(s',z)\text{d}s'\Big|(1-\tau-\frac{|s|}{s_0})^{j+1} &\leq&
\frac{s_0}{\eps}\big( 1-(1-\frac{|s|}{s_0(1-\tau) })^{j+1}\big)d_{j+1}\nonumber\\
&\leq&Cs_0d_{j+1}.\nonumber
\end{eqnarray}
This yields $ h^{-1}\|\partial_s^{-1}(b-b_0)\|_1\leq Cs_0\|b\|_1$ for
fixed $\eps>h$.\\[5pt]
Proof of $ii)$\\
By integration of inequality \eqref{ineg1}, we obtain for all $j\geq
1$
\begin{eqnarray*}
\Big|W(z)\int_0^{s} b_j(s',z)\text{d}s'\Big|
&\leq&Cs_0\frac{j!}{j\eps^j}d_j\big((1-\tau-\frac{|s|}{s_0})^{-j}-(1-\tau)^{-j}\big).\nonumber\\
&\leq&Cs_0\frac{j!}{\eps^j}d_j(1-\tau-\frac{|s|}{s_0})^{-j},\nonumber\\
\end{eqnarray*}
hence 
\begin{equation}\label{i1}
\frac{\eps^j}{j!}\Big|W(z)\int_0^{s}
b_j(s',z)\text{d}s'\Big|\Big(1-\tau-\frac{|s|}{s_0}\Big)^{j+\theta}\leq Cs_0d_j.
\end{equation}
For $j=0$ we obtain 
\begin{equation*}
\Big|W(z)\int_0^{s} b_0(s',z)\text{d}s'\Big|
\leq Cs_0\Big(\log{(1-\tau)}-\log{(1-\tau-\frac{|s|}{s_0})}\Big),\nonumber\\
\end{equation*}
then
\begin{equation}\label{i2}
\Big|W(z)\int_0^{s}
b_0(s',z)\text{d}s'\Big|\Big(1-\tau-\frac{|s|}{s_0}\Big)^{\theta}\leq Cs_0d_0.
\end{equation}
By the definition of $\|\cdot \|_{\theta}$, inequalities \eqref{i1} and
 \eqref{i2} give the result.\\[5pt]
 The proof of $iii)$  is similar, and
is left here.
\end{proof}

\begin{lemm}\label{lem3}
 There exists $C>0$ such that for all $h\in]0,\eps[$ and  $b^1,b^2\in \mathcal{H}$
\begin{equation}\label{estproduit}
\|(\partial_s^{-1}b^1)\,(\partial_s^{-1}b^2)\|_1\leq Cs_0^2\|b^1\|_1\|b^2\|_1.
\end{equation}
\end{lemm}

\begin{proof}
Write
\begin{eqnarray*}
\|(\partial_s^{-1}b^1)\,(\partial_s^{-1}b^2)\|_1&=&\sum_{j=0}^{\infty}\frac{\eps^j}{j!}\sup_{\Omega}\big|W(z)\sum_{k=0}^j\partial_s^{-1}b_k^{1}\,\partial_s^{-1}b_{j-k}^{2}(s,z)\big|\big(1-\tau-\frac{|s|}{s_0}\big)^{j+1}\\
&\leq&C\sum_{k=0}^{\infty}\sum_{j=k}^{\infty}\frac{\eps^k}{k!}\sup_{\Omega}|W(z)\partial_s^{-1}b^1_k(s,z)|\big(1-\tau-\frac{|s|}{s_0}\big)^{k+\frac12}\\
&&\quad\quad\cdot\frac{\eps^{j-k}}{(j-k)!}\sup_{\Omega}|W(z)\partial_s^{-1}b^2_{j-k}(s,z)|\big(1-\tau-\frac{|s|}{s_0}\big)^{j-k+\frac12}\\
&=C &\|\partial_s^{-1}b^1\|_{\frac12}\|\partial_s^{-1}b^2\|_{\frac12},
\end{eqnarray*}
and then by Lemma \ref{lem2} $ii)$ with $\theta=1/2$, we deduce
\begin{equation*} 
\|(\partial_s^{-1}b^1)\,(\partial_s^{-1}b^2)\|_1\leq Cs_0^2\|b^1\|_1\|b^2\|_1.
\end{equation*}
\end{proof}

\begin{prop}\label{propformel}
Let $S^0\in \mathcal{H}^0(l,B)$ be a real analytic function, and let
$a^0\in \mathcal{H}^0(l,B)$ be an analytic symbol. Then there exist
$s_0>0$, a real analytic function $S\in\mathcal{H}(s_0,l,B)$, and
an analytic symbol $a\in\mathcal{H}(s_0,l,B)$, such that 
$v=a\e^{iS/h}$ is a formal solution of equation \eqref{nlsh} with
Cauchy data $v_0=a^0\e^{iS^0/h}$.
\end{prop}

\begin{rema}
By the Cauchy formula, the function $v=a\e^{iS/h}$ satisfies for all
$k\in\N$
\begin{equation}\label{estexp1}
\sup_{|s|<s_0}\;\sup_{|\text{Im\;z}|<l/2}\;\big|\big(1-h^2\Delta\big)^{k/2}v\big|\lesssim \e^{-|z|},
\end{equation}
and
\begin{equation}\label{estexp2}
\sup_{|s|<s_0}\;\sup_{|\text{Im\;z}|<l/2}\;\big|\big(1-h^2\Delta\big)^{k/2}|v|^{p-1}v\big|\lesssim \e^{-p|z|}.
\end{equation}
This will be usefull in the sequel.
\end{rema}
\begin{proof}
The proof is based on a fixed point argument in $(\mathcal{H},\|\cdot\|_1)$.\\
Set $\varphi =\nabla S$ and differentiate the first equation of \eqref{systgeom} with respect to
the space variable, then we obtain
 
\begin{equation}\label{systgeom2}
\left\{
\begin{aligned}
&\partial_s \varphi=-2\varphi\cdot \nabla \varphi-\omega\nabla f(a_0)\\
&\partial_s a=-2\varphi\cdot\nabla
a-a\,\text{div}\varphi+ih\Delta a -\frac{i\omega a}{h}\big(f(a)-f(a_0)\big),
\end{aligned}
\right.
\end{equation}
where $\overline{a}(s,z)=\overline{a(\overline{s},\overline{z})}$ and $f(b)=(b\,\overline{b})^{\frac{p-1}{2}}$\\
Differentiate the system \eqref{systgeom2} with respect to $s$ and
obtain

\begin{equation}\label{systgeom3}
\left\{
\begin{alignedat}{2}
&\partial^2_s \varphi&=&-2\partial_s\varphi\cdot \nabla \varphi -2\varphi\cdot \nabla\partial_s \varphi -\omega\partial_s\nabla f(a_0)\\
&\partial^2_s a&=&-2\partial_s\varphi\cdot\nabla -2\varphi\cdot\nabla
  \partial_s a
-a\,\text{div}\partial_s\varphi-\partial_sa\,\text{div}\varphi+ih\Delta\partial_s a\\
&&&-\frac{i\omega\partial_s a}{h}\big(f(a)-f(a_0)\big) -\frac{i\omega a}{h}\partial_s\big(f(a)- f(a_0)\big).
\end{alignedat}
\right.
\end{equation}
Write 
\begin{equation}\label{F1}
\left\{
\begin{aligned}
&\partial_s \varphi=\partial_s^{-1} (\partial^2_s
\varphi)+\partial_s\varphi(0,\cdot),   \\
&\partial_s a=\partial_s^{-1} (\partial^2_s
a)+\partial_s a(0,\cdot),  
\end{aligned}
\right.
\end{equation}
 and
\begin{equation}\label{F2}
\left\{
\begin{aligned}
 &\varphi=\partial_s^{-2} (\partial^2_s
\varphi)+s\partial_s\varphi(0,\cdot)+\varphi(0,\cdot),   \\
&a=\partial_s^{-2} (\partial^2_s
a)+s\partial_s a(0,\cdot)+a(0,\cdot).  
\end{aligned}
\right.
\end{equation}
Now introduce the new unknown function $u=(u_1,u_2)=(\partial_s^2 \varphi,
\partial_s^2 a)$. Hence we are lead to solving  a system of the form 
\begin{equation}\label{ptfixe}
u=F(s,u).
\end{equation}
We will show that for $0<s_0<1$ small enough $F$ is a contraction in a
ball in $(\mathcal{H},\|\cdot\|_1)$. Let $R>0$ be such that 

\begin{equation*}
\|\varphi(0,\cdot)\|_0,\|\partial_s\varphi(0,\cdot)\|_0,\|\nabla\varphi(0,\cdot)\|_0,
\|\nabla\partial_s\varphi(0,\cdot)\|_0,\|\Delta\partial_s\varphi(0,\cdot)\|_0
\leq R,
\end{equation*}
and
\begin{eqnarray*}
&\|a(0,\cdot)\|_0,\|\partial_sa(0,\cdot)\|_0,\|\nabla a(0,\cdot)\|_0,
\|\nabla\partial_s a(0,\cdot)\|_0,\|\Delta\partial_s a(0,\cdot)\|_0,\\
&\|{(a-a_0)(0,\cdot)}/{h}\|_0,\|{\partial_s(a-a_0)(0,\cdot)}/{h}\|_0\leq R.
\end{eqnarray*}
$\bullet$ Write 
\begin{eqnarray*}
  \partial_s \varphi\nabla\varphi&=&
 \big( \partial_s^{-1}  \partial^2_s \varphi+
\partial_s\varphi(0,\cdot)    \big)\cdot\\
&&\quad \big(  \partial_s^{-1}(
\partial_s^{-1}\nabla)(\partial^2_s \varphi)+s\nabla
\partial_s\varphi(0,\cdot)+ \nabla\varphi(0,\cdot) \big) 
\end{eqnarray*}
Then by \eqref{estproduit} and \eqref{estnorm}
\begin{eqnarray*}
\| \partial_s \varphi\nabla\varphi\|_1&\lesssim&s_0^2\| \partial^2_s
\varphi\|_1\| \partial_s^{-1}\nabla(\partial^2_s \varphi)  \|_1+R\|
\partial_s^{-1}\partial^2_s \varphi  \|_1\\
&&+R\| \partial_s^{-1} (\partial_s^{-1}\nabla)(\partial^2_s \varphi)  \|_1+R^2,
\end{eqnarray*}
and by \eqref{inegalite1} and \eqref{inegalite2}
\begin{eqnarray}
\| \partial_s \varphi\nabla\varphi\|_1&\lesssim&s_0^3\| \partial^2_s
\varphi\|_1^2+s_0R\| \partial^2_s\varphi\|_1+R^2\nonumber\\
&\lesssim&s_0^2\| \partial^2_s
\varphi\|_1^2+R^2\label{£1}.
\end{eqnarray}
$\bullet$ Similarly we obtain 
\begin{equation}\label{£2}
\| \varphi\nabla \partial_s\varphi\|_1\lesssim s_0^2\| \partial^2_s\varphi\|_1^2+R^2,
\end{equation}
and 
\begin{equation}\label{£3}
\| \partial_s \varphi\nabla a\|_1, \| \varphi\nabla \partial_s a\|_1,
\|a\,\text{div}\, \partial_s \varphi\|_1,\|\partial_sa\,\text{div}\, \varphi\|_1
\lesssim s_0^2\| \partial^2_s
\varphi\|_1\| \partial^2_s a\|_1+R^2.
\end{equation}
$\bullet$ We have
 \begin{equation}\label{z0}
\nabla \partial_sf(a_0)=\nabla \partial_s^{-1} \partial_s^{2}f(a_0)+\nabla \partial_sf(a_0)(0,\cdot).
\end{equation}
By the Leibniz rule and \eqref{estnorm}
\begin{equation*}
\| \partial^2_s f(a_0)\|_1\lesssim \| \partial^2_s a_0\|_1 \| a_0\|_0^{{p-2}}
+ \| (\partial_s a_0)^2\|_1 \| a_0\|_0^{{p-3}}.
\end{equation*}
From \eqref{F2} and \eqref{inegalite3}
 \begin{equation}\label{z1}
\| a_0\|_0\lesssim \| \partial_s^{-2} \partial^2_s a_0\|_0+R\lesssim s_0\| \partial^2_s a_0\|_1+R.
\end{equation}
Now use \eqref{F1}, \eqref{estproduit} and \eqref{inegalite2}
 \begin{equation}\label{z2}
\| (\partial_s  a_0)^2\|_1\lesssim s_0^2\| \partial^2_s a_0\|_1^2+R^2.
\end{equation}
Finally, from \eqref{z0}, \eqref{z1} and \eqref{z2} we deduce 
\begin{equation}\label{£4}
\|\nabla \partial_sf(a_0) \|_1\lesssim s_0^{\frac{p-1}2} \|
\partial^2_s a_0\|_1^{\frac{p-1}2}+R^{\frac{p-1}2}
\lesssim s_0^{\frac{p-1}2} \| \partial^2_s a\|_1^{\frac{p-1}2}+R^{\frac{p-1}2}.
\end{equation}
$\bullet$ Write
\begin{equation*}
h\Delta \partial_s a=h \partial_s^{-1}\Delta  \partial^2_s a+h\Delta \partial_s a(0,\cdot),
\end{equation*}
therefore by \eqref{inegalite1} 
\begin{equation}\label{£5}
\|h\Delta \partial_s a\|_1 \lesssim s_0 \| \partial^2_s a\|_1+R.
\end{equation}
$\bullet$ We now estimate the term $\frac{\partial_s
  a}{h}\big(f(a)-f(a_0)\big)$. Observe that 
\begin{eqnarray*}
f(a)-f(a_0)&=&(a\overline{a})^{\frac{p-1}2}-(a_0\overline{a_0})^{\frac{p-1}2}\\
&=&(a\overline{a}-a_0\overline{a_0} )\big(  (a\overline{a})^{\frac{p-3}2} +\cdots+  (a_0\overline{a_0})^{\frac{p-3}2}   \big),
\end{eqnarray*}
and 
\begin{equation*}
a\overline{a}-a_0\overline{a_0} =(a-a_0)\overline{a}+(\overline{a}-\overline{a_0})a_0.
\end{equation*}
Then by \eqref{estnorm}
\begin{equation}\label{*1}
\|\frac{\partial_s a}{h}\big(f(a)-f(a_0)\big)\|_1\lesssim \|\partial_s a\,\frac{a-a_0}{h}\|_1\|a\|_0^{p-2}.
\end{equation}
Use \eqref{F1},  \eqref{F2} to write 
\begin{eqnarray*}
\partial_s a\,\frac{a-a_0}{h}&=&\Big( \partial_s^{-1}(\partial_s^2a)+
\partial_sa(0,\cdot)      \Big)\\
&&\quad \quad \Big( \partial_s^{-1}  \partial_s^{-1} \frac{\partial_s^{2}(a-a_0)}{h} 
+s\frac{\partial_s(a-a_0)(0,\cdot)}{h}+ \frac{(a-a_0)(0,\cdot)}{h}  \Big),
\end{eqnarray*}
then by \eqref{estproduit}, \eqref{estnorm} and \eqref{inegalite2}

\begin{equation}\label{*2}
\|\partial_s a\,\frac{a-a_0}{h}\|_1\lesssim \big(s_0\| \partial_s^{2}a\|_1+Rs_0 \big)
 \big( \|   \partial_s^{-1}   \frac{\partial_s^{2}(a-a_0)}{h}\|_1  +R     \big).
\end{equation}
Moreover from \eqref{inegalite1} we have 
\begin{equation}\label{*3}
 \|   \partial_s^{-1}   \frac{\partial_s^{2}(a-a_0)}{h}\|_1 \lesssim s_0\| \partial_s^{2}a\|_1.
\end{equation}
Therefore inequalities \eqref{*1}, \eqref{*2} and  \eqref{*3} yield
\begin{equation}\label{£6}
 \|\frac{\partial_s a}{h}\big(f(a)-f(a_0)\big)\|_1\lesssim s_0^p\|\partial_s^2a\|_1^p+R^p.
\end{equation}
$\bullet$ Similar arguments are used to show that 
\begin{equation}\label{£7}
 \|\frac{a}{h}\partial_s\big(f(a)- f(a_0)\big)\|_1\lesssim s_0^p\|\partial_s^2a\|_1^p+R^p.
\end{equation}

\noindent Inequalities \eqref{£1}, \eqref{£2}, \eqref{£3}, \eqref{£4}, \eqref{£5}, \eqref{£6}
and \eqref{£7} show that, if $s_0>0$ is small enough, there exists
$R_1>R$ such that $F$ maps the ball of radius $R_1$ \big(in
($\mathcal{H},\|\cdot\|_1$)\big) into itself.\\[5pt]
With analogous arguments, we can show that $F$ is a contraction in
($\mathcal{H},\|\cdot\|_1$).\\
Hence by the fixed  point theorem, there exists a unique $u=(\partial^2_s\varphi,\partial^2_sa)\in
\mathcal{H}\times\mathcal{H} $ which satisfies \eqref{ptfixe}.\\
Let $(\varphi^0,a^0)\in\mathcal{H}\times\mathcal{H}$, and consider the 
couple $(\partial_s\varphi(0,\cdot),\partial_s
a(0,\cdot))\in\mathcal{H}\times\mathcal{H}$ which solves the system
\eqref{systgeom2} at $s=0$. Let $u$ be the solution of \eqref{ptfixe}
with these initial conditions. Then with the formula \eqref{F2} we
recover the couple $(\varphi,a)$ which is a solution of \eqref{systgeom2}. Moreover, \eqref{F2} shows
that $(\varphi,a)\in\mathcal{H}\times\mathcal{H} $.\\[5pt]
Let $S^0\in \mathcal{H}^0$ and take $\varphi^0=\nabla S^0$. The
function $\varphi$ (with Cauchy condition
$\varphi(0,\cdot)=\varphi^0$) is irrotational, as it satisfies the
equation
$$\partial_s \varphi=-2\varphi\cdot \nabla \varphi-\omega\nabla
f(a_0). $$
Therefore there exists $S$ so that $\varphi=\nabla S$ and which is
solution of 
\begin{equation*}
\nabla\big(\partial_s S+(\nabla S)^2+\omega f(a_0)\big)=0.
\end{equation*}
Moreover, it is possible to choose $S$  such that 
\begin{equation*}
\partial_s S+(\nabla S)^2+\omega f(a_0)=0.
\end{equation*}
Now the formula 
\begin{equation}\label{F3}
S(s,z)=\int_{0}^s\partial_sS(\sigma, z)\text{d}\sigma+S^0(z)
=-\int_{0}^s\big(\varphi\cdot \varphi+\omega f(a_0)\big)(\sigma, z)\text{d}\sigma+S^0(z),
\end{equation}
shows that  $S\in \mathcal{H}$.\\[5pt]
Finally, we have shown the existence of a solution $(S,a)\in
\mathcal{H}\times\mathcal{H}$ of the system 
\begin{equation*}
\left\{
\begin{aligned}
&\partial_s S+(\nabla S)^2+\omega|a_0|^{p-1}=0,\\
&\partial_s a+2\nabla S\cdot\nabla
a+a\Delta S-ih\Delta a +\frac{i\omega a}{h}(|a|^{p-1}-|a_0|^{p-1})=0,\\
&S(0,z)=S^0(z)\in \mathcal{H}^0,\;a(0,z,h)=a^0(z,h)\in \mathcal{H}^0.
\end{aligned}
\right.
\end{equation*}
With a Gronwall inequality, it is straightforward to check that $S$ is
real analytic.
\end{proof}
\begin{rema}
The inequality $\|\partial_s^{-1}b\|_0\leq Cs_0\|b\|_1$ fails, and that is
the reason why we have to differentiate the system \eqref{systgeom2}
with respect to the time variable, before applying the contraction method.
\end{rema}


\subsection{The general case of an analytic manifold
  $(M^d,g)$}\label{SectionAna}
$~$\\
Let $(M^d,g)$ an analytic Riemannian manifold of dimension $d$. We
assume moreover that $g$ is analytic. Let $m\in
M^d$. Then there exist a  neighbourhood $\mathcal{U}\subset M^d$ of
$m$, a neighbourhood $\mathcal{V}\subset \R^d$ of $0$, and an
homeomorphism 
\begin{equation}\label{defkappa}
\kappa\colon \mathcal{U}\longrightarrow \mathcal{V}.
\end{equation}
In the chart $(\mathcal{U},\kappa)$ the metric $g$ can be written
$$g=\sum_{1\leq j,k\leq d}g_{jk}(x)\text{d}x_j\text{d}x_k,$$
where $G=(g_{jk})$ is a positive symmetric matrix and analytic in
$\mathcal{V} $.\\
In these coordinates, we have the explicit formula for the
Laplace-Beltrami operator
\begin{eqnarray*}
\Delta_g =\Delta_g(x)&=&\frac{1}{\sqrt{\text{det}G}}\text{div} \big( \sqrt{\text{det}G}\;G^{-1}\nabla\cdot \big)\\
&=&\frac{1}{\sqrt{\text{det}G}}\sum_{1\leq j,k\leq
  d}\frac{\partial}{\partial x_j}\big(\sqrt{\text{det}G}\;g^{jk}\frac{\partial}{\partial x_k}\big),
\end{eqnarray*}
where $(g^{jk})=G^{-1}$. Every function involved in the former
expression is analytic.\\
We now make the rescaling \eqref{chgt}. The function
$$v(s,z,h)=\h^{-\gamma}u(\h^{\alpha}s,\h z), $$
satisfies
\begin{equation}\label{Mnls}
i h\partial_t v(s,z)+h^2\Delta(\h z) v(s,z)   =  \omega|v|^{p-1} v(s,z),\quad
(s,z)\in\R\times \h^{-1}\mathcal{V}.
\end{equation}
We now adapt the analysis of  Section \ref{SectionWKB} to the
equation \eqref{Mnls}, in $\h^{-1}\mathcal{V}$ instead of
$\R^d$.\\
Let $r>0$ such that 
\begin{equation}\label{defr}
\mathcal{B}(0,2r)\subset\mathcal{V}.
\end{equation}
Notice that on the set $\big\{(|\h z|<r)\cap(|\text{Im}\;z|<l)\big\}$,
the coefficients of $\Delta_g$ are uniformly bounded with respect to
$\h$, as well as their derivatives.\\
Here again, we want to find a formal solution of \eqref{Mnls} of the form
\begin{equation*}
v(s,z,h)= a(s,z,h)\e^{iS(s,z)/h}=\big(\sum_{j\geq 0} a_j(s,z)h^j\big)\e^{iS(s,z)/h}.
\end{equation*}
Therefore $(S,a)$ has to satisfy  the system
\begin{equation}\label{Msystgeom}
\left\{
\begin{aligned}
&\partial_s S+(\nabla_g S)^2+\omega|a_0|^{p-1}=0,\\
&\partial_s\, a+2\nabla_g S\cdot\nabla_g\,
a+a\,\Delta_g S-ih\Delta_g \,a +\frac{i\omega a}{h}(|a|^{p-1}-|a_0|^{p-1})=0,\\
&S(0,z)=S^0(z),\;a(0,z,h)=a^0(z,h),
\end{aligned}
\right.
\end{equation}
with $\nabla_g=\nabla_g(\h z)$, $\Delta_g=\Delta_g(\h z)$ and where $v(0,z,h)=a^0(z,h)\e^{iS^0(z)/h}$.\\
For $\h>0$ small enough, denote by 
\begin{equation*}
\mathcal{D}_{\h}= \{|s|<s_0\}\times
 \big\{(|\h z|<r)\cap(|\text{Im}\;z|<l)\big\},
\end{equation*}
and by  $\mathcal{H}_{\h}=\mathcal{H}_{\h}(s_0,l,r,B)$ the space of the analytic symbols
 $b(s,z,h)=\sum_{j\geq 0}b_j(s,z)h^j$ (see Definition \ref{def1}) satisfying
\begin{equation*}
|W(z)\,b_j(s,z)|\leq AB^jj!\;\;\text{on}\;\;\mathcal{D}_{\h},\;\forall j\geq 0.
\end{equation*}
Define also $\mathcal{H}^0_{\h}=\mathcal{H}_{\h}(0,l,r,B)$ the space of
the initial conditions.

\noindent Let $\eps<1/B$. For $0\leq \theta\leq 1 $, we  endow $\mathcal{H}_{\h}(s_0,l,r,B)$ with the norms 
\begin{equation*}
\|b\|_{\theta,\h}=\sum_{j\geq 0}\frac{\eps^j}{j!}\sup_{0<\tau<1}\;
\sup_{|s|<s_0(1-\tau)}\;\sup_{\Gamma_{\tau}}|W(z)\,b_j(s,z)|\Big(1-\tau-\frac{|s|}{s_0}\Big)^{j+\theta},
\end{equation*}
where $\Gamma_{\tau}=\big\{(|\h z|<r\tau)\cap(|\text{Im}\;z|<l\tau)\big\}$.\\[5pt]
Now it is straightforward to check that the results of Lemma
\ref{lem1}, Lemma \ref{lem2} and  Lemma \ref{lem3} hold when
$\|\cdot\|_{\theta}$ is replaced with $\|\cdot\|_{\theta,\h}$ and that
the constants involved in the estimates do not depend on $\h$. Notice
that  the boundedness of $\mathcal{D}_{\h}$ with respect to the
variable $|\h z|$ is dealt with in exaclty the same way as was done
with the boundedness with respect to $|\text{Im}\;z|$. This
yields the following analog of Proposition \ref{propformel}

\begin{prop}\label{Mpropformel}
Let $S^0\in \mathcal{H}_{\h}^0(l,r,B)$ be a real analytic function, and let
$a^0\in \mathcal{H}_{\h}^0(l,r,B)$ be an analytic symbol. Then there exist
$s_1>0$ independent of $\h$, a real analytic function $S\in\mathcal{H}_{\h}(s_1,l,r,B)$, and
an analytic symbol $a\in\mathcal{H}_{\h}(s_1,l,r,B)$, such that 
$v=a\e^{iS/h}$ is a formal solution of equation \eqref{nlsh} with
Cauchy data $v_0=a^0\e^{iS^0/h}$.
\end{prop}

\begin{rema}
Proposition \ref{propformel} is contained in Proposition
\ref{Mpropformel}: In the case $(M^d,g)=(\R^d,\text{can})$, $r=+\infty$ and 
$\mathcal{H}_{\h}(s_1,l,r,B)=\mathcal{H}(s_1,l,B)$.
\end{rema}
\noindent We  are now able to  construct  an  approximate solution of
the problem \eqref{Mnls}.\\[5pt]
 Let $c_0$ such that $c_0/h=:n\in\N$.
Define 

\begin{equation*}
a^{(n)}(s,z,h)=\sum_{j\leq n}a_j(s,z)h^j,
\end{equation*}
and 
\begin{equation}\label{defvapp}
v_{\text{app}}(s,z,h)=a^{(n)}(s,z,h)\e^{iS(s,z)/h}.
\end{equation}
where the $a_j$'s and $S$ are given by Proposition
\ref{Mpropformel}. The choice of the initial condition
$v_{\text{app}}(0,z,h)$ will be made in Section \ref{SectionInstab}.\\[5pt]
We now show that if $c_0$ is small enough, $v_{\text{app}}$ is a good
approximation to the problem \eqref{nlsh}.
\begin{prop}\label{corapproxsc}
Let $s_1>0$ be given by Proposition \ref{Mpropformel}. If $c_0\ll1$, there exists
$\delta_1>0$ such that the function $v_{\text{app}}$ defined by \eqref{defvapp} satisfies
\begin{equation}
i h\partial_sv_{\text{app}} +h^2\Delta v_{\text{app}}
=  \omega (v_{\text{app}} \overline{v_{\text{app}}})^{\frac{p-1}{2}} v_{\text{app}}
+\e^{-\delta_1/h}g,
\end{equation}
 with $\overline{v_{\text{app}}}=\overline{v_{\text{app}}(s,\overline{z})}$ and where $g$ is an analytic function on $\{|s|<s_1\}\times \big\{(|\h
 z|<r)\cap(|\text{Im}\;z|<l)\big\}$ such that for all $k\in \N$, there exists $C_k>0$
independent of $h$ so that 
\begin{equation}\label{ee1}
\sup_{|s|<s_1}\;\sup_{|\text{Im}\;z|<l/2}\|(1-h^2\Delta)^{k/2}g(s,\cdot+i\text{Im}\;z)
\|_{L^2(\mathcal{B}(0,r/\h))}\leq C_k.
\end{equation}
\end{prop}
\noindent Here we have used the convention that $\mathcal{B}(0,r/\h)=\R^d$ if $r=+\infty$. 
\begin{proof} Denote by  $f(b)=(b\,\overline{b})^{\frac{p-1}{2}}$, with $\overline{b}=\overline{b(s,\overline{z})}$. The function  $v_{\text{app}}$ satisfies the equation
\begin{eqnarray}\label{eqansatz}
\begin{array}{l}
i h\partial_sv_{\text{app}}+h^2\Delta v_{\text{app}}
- \omega f(v_{\text{app}})
v_{\text{app}}\hspace{30pt}\\[5pt]
\;=-a^{(n)}\big(     \partial_s S+(\nabla S)^2+\omega f(a_0)\big)\e^{iS/h}\\[5pt]
\hspace{30pt}+ih\Big(\partial_s a^{(n)}+2\nabla S\cdot\nabla
a^{(n)}+a^{(n)}\Delta S-ih\Delta a^{(n)} \\[5pt]
\hspace{60pt}+\displaystyle\frac{i\omega
  a^{(n)}}{h}\big(  f(a^{(n)})-f(a_0) \big)      \Big)\e^{iS/h}.
\end{array}
\end{eqnarray}
For $m=n,n+1$ write the expansion in $h$
\begin{equation}\label{bqq}
\frac{i\omega
  a^{(m)}}{h}\big(  f(a^{(m)})-f(a_0)  \big) :=\sum_{j=0}^{pm-1}b_{j,m}h^{j}.   
\end{equation}
By construction the following system is satisfied

\begin{equation}\label{systN}
\left\{
\begin{aligned}
&\partial_s S+(\nabla S)^2+\omega(a_0\overline{a_0})^{\frac{p-1}{2}}=0,\\
&\partial_s a^{(n)}+2\nabla S\cdot\nabla
a^{(n)}+a^{(n)}\Delta S-ih\Delta a^{(n-1)}+\sum_{j=0}^{n}b_{j,n+1}h^{j}=0.
\end{aligned}
\right.
\end{equation}
Notice that 
\begin{equation}\label{bq}
b_{j,n}=b_{j,n+1}\quad \text{for all }\quad j\leq n-1.
\end{equation}

\noindent Therefore by \eqref{bq} and \eqref{systN}, \eqref{eqansatz} rewrites

\begin{eqnarray}
\lefteqn{ i h\partial_sv_{\text{app}}+h^2\Delta v_{\text{app}}
- \omega ( v_{\text{app}} \overline{v_{\text{app}}})^{\frac{p-1}{2}}
v_{\text{app}}}\nonumber\\
&=&ih\big( -ih^{n+1}\Delta a_{n}- \sum_{j=0}^{n}b_{j,n+1}h^{j} +
\sum_{j=0}^{pn-1}b_{j,n}h^{j} \big)\e^{iS/h}\nonumber\\
\displaystyle&=& \big( h^{n+2}\Delta a_{n}- i h^{n+1} b_{n,n+1} +
ih\sum_{j=n}^{pn-1}b_{j,n}h^{j}  \big)\e^{iS/h}.\label{eqansatz2}
\end{eqnarray}
We now estimate each term of the r.h.s. of \eqref{eqansatz2}. By \eqref{bqq} we have 
\begin{equation*}
hb_{j,n}=i\omega\Big( \sum_{i_1+\cdots +i_p=j}\widetilde{a_{i_1}}\cdots\widetilde{a_{i_p}}-(a_0\overline{a_0})^{\frac{p-1}{2}}a_j\Big),
\end{equation*}
with $\widetilde{a_{i_k}}=a_{i_k}$ or
$\widetilde{a_{i_k}}=\overline{a_{i_k}}$.\\
Now by \eqref{expanalb}, $|a_{i_k}|\lesssim B^{i_k}(i_k)!\,\e^{-|z|}$,
  thus 
\begin{equation}\label{hb}
h|b_{j,n}|\lesssim B^j\Big(  \sum_{i_1+\cdots +i_p=j}(i_1)!\,
 \cdots(i_p)!\,  +j!\, \Big)\e^{-p|z|}\lesssim  B^jj!\,\e^{-p|z|},
\end{equation}
and by the Stirling formula,
\begin{equation*}
 (pn)!\,\lesssim n^{1/2}\big(\frac{pn}{\e}\big)^{\textstyle{pn}},
\end{equation*}
we deduce from \eqref{hb}
\begin{eqnarray}
|h\sum_{j=n}^{pn-1}b_{j,n}h^j|&\lesssim&
 \big(\sum_{j=n}^{pn-1}B^jj!\,h^j\big)\e^{-p|z|}\nonumber\\
&\leq& (p-1)n(Bh)^{pn}(pn)!\, \e^{-p|z|}\nonumber\\[4pt]
&\lesssim& h^{-\frac32}\big(\frac{Bc_0p}{\e}\big)^{\displaystyle{\frac{c_0p}{h}}}\e^{-p|z|},\nonumber
\end{eqnarray}
as we have $n=c_0/h$. Now choose $c_0<e/(Bp)$, then there exists $\delta>0$ such that 
\begin{equation*}
|h\sum_{j=n}^{pn-1}b_{j,n}h^j|\lesssim \sum_{j=n}^{pn-1}h|b_{j,n}h^j|\lesssim \e^{-\delta/h}\e^{-p|z|}.
\end{equation*}
Similarly, for some $\delta>0$
\begin{eqnarray}
|h^{n+1}b_{n,n+1}|&\lesssim& \e^{-\delta/h}\e^{-|z|},\nonumber\\
|h^{n+2}\Delta a_n|&\lesssim& \e^{-\delta/h}\e^{-|z|}.\nonumber
\end{eqnarray}
Finally use that the function $\phi:(\text{Re\,}z,\text{Im\,}z)\longmapsto
\e^{-|\text{Re\,}z+i\text{Im\,}z|}$ satisfies
\begin{equation*}
\sup_{|\text{Im\,}z|<l}\|\phi(\cdot,\text{Im\,}z)\|_{L^2(\mathcal{B}(0,r/\h))}\lesssim 1.
\end{equation*}
We have therefore proved the estimate \eqref{ee1} for $k=0$.\\
To treat the case $k\geq 0$, use the Cauchy formula to obtain 
\begin{equation*}
\sup_{|s|<s_1}\;\sup_{|\text{Im\;z}|<l/2}\big| (1-h^2\Delta)^{k/2}a_j\big|\lesssim
\sup_{|s|<s_1}\;\sup_{|\text{Im\;z}|<l}|a_j|\lesssim B^jj!\,\e^{-|z|},
\end{equation*}
and 
\begin{equation*}
\sup_{|s|<s_1}\;\sup_{|\text{Im\;z}|<l/2}\big|
(1-h^2\Delta)^{k/2}\e^{iS/h}\big|\lesssim 1,
\end{equation*}
and we can easily adapt the previous computations. 
\end{proof}

\section{Validity of the Ansatz }\label{SectionAnsatz}

\begin{prop}\label{propapprox}
Let  $v_{\text{app}}$ be the function  defined by
\eqref{defvapp}. Let $v$ be the  solution of

\begin{equation}\label{eqv}
\left\{
\begin{aligned}
& i h\partial_sv +h^2\Delta v   = \omega | v|^{p-1} v,\quad
(s,z)\in\R^{1+d},\\
&v(0,z)=v_{\text{app}}(0,z).
\end{aligned}
\right.
\end{equation}
Then there exist $s_2>0$ and $\delta_2>0$ such that for all $k\in\N$
\begin{equation*}
\sup_{0<s<s_2}\;\|
(1-h^2\Delta)^{k/2}(v-v_{\text{app}})(s)\|_{L^2(\mathcal{B}(0,r/\h))}\leq C_k \e^{-\delta_2/h},
\end{equation*}
with $C_k>0$.
\end{prop}

\begin{proof} 
It is given in \cite{PGX}, but we reproduce it in the appendix.
\end{proof}

\noindent We are now able to define the Ansatz to the equation
\eqref{nls}.\\[5pt]
In the case $(M^d,g)=(\R^d,\text{can})$, we consider the function
$u_{\text{app}}$ given by \eqref{defvapp} and define 
\begin{equation}\label{defuapp22}
u_{\text{app}}(t,x)=\h^{\gamma}v_{\text{app}}(\h^{-\alpha}t,\h^{-1}x),
\end{equation} 
where $\gamma$ and $\alpha$ satisfy the relations
\eqref{parameters} and $h=\h^{\beta}$. The initial condition will be
given in the next section.\\
From Proposition \ref{propapprox} we deduce

\begin{coro}\label{corapprox}
(The case $(M^d,g)=(\R^d,\text{can})$) Let $s_2$ be given by Proposition \ref{propapprox}, let $u_{\text{app}}$ be given by \eqref{defuapp22}, and let $u$ be the
solution of
\begin{equation*}
\left\{
\begin{aligned}
&i \partial_t u+\Delta u   = \omega |u|^{p-1} u,\quad (t,x)\in \R^{1+d},\\
&u(0,x)=u_{\text{app}}(0,x).
\end{aligned}
\right.
\end{equation*}
 Then for all $k\in\N$
\begin{equation*}
\|u-u_{\text{app}}\|_{L^{\infty}([0,\h^{\alpha}s_2];H^k(\R^d))}\longrightarrow
0,
\end{equation*}
when $h\longrightarrow 0$.
\end{coro}

\noindent In the general case of an analytic manifold $(M^d,g)$, we
have to construct an approximate solution supported in
$\mathcal{B}(0,r)\subset \mathcal{U}$.\\
Let $\chi\in \mathcal{C}^{\infty}_0(\R)$, $\chi\geq0$, such that
\begin{equation}\label{defchi}
\chi(\xi)=\left\{\begin{array}{l} 
1\quad \text{for} \;|\xi|\leq r/2, \\[5pt]   
0 \quad \text{for} \;|\xi|\geq r.
\end{array} \right.
\end{equation}
 Let $0<\eta<1$, let $v_{\text{app}}$ be given by \eqref{defvapp} and
consider
\begin{equation}\label{defuapp}
\underline{u_{\text{app}}}(t,x)=\h^{\gamma}\chi(\h^{-\eta}|x|)v_{\text{app}}(\h^{-\alpha}t,\h^{-1}x),
\end{equation} 
where $\gamma$ and $\alpha$ are given by the relations
\eqref{parameters}, and $h=\h^{\beta}$.\\
We have
\begin{equation*}
\text{supp\;}\underline{u_{\text{app}}}\subset \{(t,x)\in\R^{1+d}\;:\;\;|x|\leq r\h^{\eta}\},
\end{equation*} 
which concentrates in $x=0$.\\
Hence if $\h$ is small enough, $u_{\text{app}}$ is supported in
$\mathcal{V}$, and we can transport this function  to $\mathcal{U}$ by
the chart $\kappa$ (see \eqref{defkappa}).We therefore define the
approximate solution $u^M_{\text{app}}$ of \eqref{nls} by
\begin{equation}\label{Muapp}
u^M_{\text{app}}=\underline{u_{\text{app}}}\circ \kappa.
\end{equation}
 In the following we write $u^M_{\text{app}}=u_{\text{app}}$.\\
Then, as ${u}_{\text{app}}$ is compactly supported, it can not be analytic. We
now consider all the  functions only with real variables.\\
Up to now, we did not use the rate of decrease of the weight
$W^{-1}$ introduced in \eqref{expanalb}, but it is needed now because
of the truncation. However, because of the error $\e^{-c/\h}$ induced
from this cutoff, we obtain the following weaker result

\begin{coro}\label{Mcorapprox}
(The general case) Let $s_2$ be given by Proposition \ref{propapprox}, let $u_{\text{app}}$ be given by \eqref{defuapp}, and let $u$ be the
solution of
\begin{equation}\label{nlsapp}
\left\{
\begin{aligned}
&i \partial_t u+\Delta u   = \omega |u|^{p-1} u,\quad (t,x)\in \R
  \times M^d,\\
&u(0,x)=u_{\text{app}}(0,x).
\end{aligned}
\right.
\end{equation}
Let $\kappa\ge 0$ such that $\beta+\eta-\kappa< 1$. Then for all $k\in\N$
\begin{equation*}
\|u-u_{\text{app}}\|_{L^{\infty}([0,\h^{\alpha+\kappa}s_2];H^k(M^d))}\longrightarrow
0,
\end{equation*}
when $h\longrightarrow 0$.
\end{coro}

\begin{proof}
Let $k>d/2$ an integer, and set 
\begin{equation*}
\|f\|_{H^k_{\h}}=\|\big(1-{\h}^{2(\beta+1)}\Delta\big)^{k/2}f\|_{L^2(M^d)}.
\end{equation*}
With the Leibniz rule and interpolation we check that for all
$f\in H^{k}(M^d)$ and $g\in W^{k,\infty}(M^d)$
\begin{equation}\label{produithh}
\|f\,g\|_{H^k_{\h}}\lesssim
\|f\|_{H^k_{\h}}\|g\|_{L^{\infty}(M^d)}+\|f\|_{L^2(M^d)}
\|\big(1-{\h}^{2(\beta+1)}\Delta\big)^{k/2}g\|_{L^{\infty}(M^d)}.
\end{equation}
Moreover, as $k>d/2$, for all $f_1,f_2\in {H^k(M^d)}$
\begin{equation}\label{produithh2}
\|f_1\,f_2\|_{H^k_{\h}}\lesssim   \h^{-(\beta+1)k}\|f_1\|_{H^k_{\h}}\|f_2\|_{H^k_{\h}}
\end{equation}
The function $u_{\text{app}}$ satisfies
\begin{equation*}
i \partial_t u_{\text{app}}+\Delta u_{\text{app}}   = 
\omega |u_{\text{app}}|^{p-1} u_{\text{app}}+\e^{-{c}/{\h^{1-\eta}}}q,
\end{equation*}
with
 \begin{equation*}
\|q\|_{H^k_{\h}}\lesssim 1.
\end{equation*}
Let $u$ be the solution of \eqref{nlsapp} and define
$w=u-u_{\text{app}}$. Then $w$ satisfies
\begin{equation}\label{nlsapp2}
\left\{
\begin{aligned}
&i \partial_t w+\Delta w   = \omega\big( |w+u_{\text{app}}|^{p-1}
  (w+u_{\text{app}})- |u_{\text{app}}|^{p-1} u_{\text{app}} \big)+\e^{-{c}/{\h^{1-\eta}}}q\\
&w(0,x)=0.
\end{aligned}
\right.
\end{equation}
We expand the r.h.s. of \eqref{nlsapp2}, apply the operator
$\big(1-{\h}^{2(\beta+1)}\Delta\big)^{k/2}$ to the equation, and take
the $L^2$- scalar product with
$\big(1-{\h}^{2(\beta+1)}\Delta\big)^{k/2}w$. Then we obtain
\begin{equation}\label{nlsapp5}
\frac{\text{d}}{\text{d}t}\|w\|_{H^k_h}\lesssim \sum_{j=1}^p
 \|w^j\,u^{p-j}_{\text{app}}\|_{H^k_h}+\e^{-{c}/{\h^{1-\eta}}}.
\end{equation}
We now have to estimate the terms
$\|w^j\,u^{p-j}_{\text{app}}\|_{H^k_h},$
for $1\leq j\leq p$. From \eqref{produithh} we deduce 
\begin{eqnarray}\label{nlsapp3}
\|w^j\,u^{p-j}_{\text{app}}\|_{H^k_h}&\lesssim&
\|w^j\|_{H^k_h}\|u^{p-j}_{\text{app}}\|_{L^{\infty}(M^d)}\nonumber\\
&&\quad+\|w^j\|_{L^2(M^d)}\|\big(1-{\h}^{2(\beta+1)}\Delta\big)^{k/2}u^{p-j}_{\text{app}}\|_{L^{\infty}(M^d)}.
\end{eqnarray}
By \eqref{produithh2}, and as we have
\begin{equation}\label{nlsapp4}
\|u^{p-j}_{\text{app}}\|_{L^{\infty}(M^d)}\lesssim  \h^{\gamma(p-j)},\quad
\|\big(1-{\h}^{2(\beta+1)}\Delta\big)^{k/2}u^{p-j}_{\text{app}}\|_{L^{\infty}(M^d)}\lesssim \h^{\gamma(p-j)},
\end{equation}
thus inequality \eqref{nlsapp3} yields
\begin{equation*}
\|w^j\,u^{p-j}_{\text{app}}\|_{H^k_h}\lesssim  \h^{\gamma(p-j)} \h^{-(\beta+1)(j-1)k}\|w\|^{j}_{H^k_h}.
\end{equation*}
Therefore, from \eqref{nlsapp5} we have
\begin{equation*}
\frac{\text{d}}{\text{d}t}\|w\|_{H^k_h}\lesssim\h^{\gamma(p-1)}\|w\|_{H^k_h}
+  \h^{-(\beta+1)(p-1)k}\|w\|^{p}_{H^k_h}+\e^{-{c}/{\h^{1-\eta}}}.
\end{equation*}
Observe that $\|w(0)\|_{H^k_h}=0$. Now, for times $t$ so that 
\begin{equation}\label{bootstrap}
\|w\|_{H^k_h}\lesssim \h^{\gamma+(\beta+1)k},
\end{equation}
we can remove the nonlinear term in \eqref{nlsapp3}, and by the
Gronwall Lemma,
 \begin{equation}\label{bootstrap2}
\|w\|_{H^k_h}\lesssim \e^{-{c}/{\h^{1-\eta}}}\e^{C\h^{\gamma(p-1)}t}.
\end{equation}
By \eqref{parameters}, $\alpha=2+\beta$ and $\gamma(p-1)=-2(\beta+1)$,
thus for all $0\leq t\leq s_2\h^{\alpha+\kappa}$,
$$ \h^{\gamma(p-1)}t\leq s_2\h^{-\beta+\kappa},$$
and if $\beta+\eta-\kappa<1$, the r.h.s. in \eqref{bootstrap2} tends to
$0$. Then the inequality \eqref{bootstrap} is satisfied for all
$0\leq t\leq s_2\h^{\alpha+\kappa}$, and 
with a continuity argument, we infer that \eqref{bootstrap2} holds for
$0\leq t\leq s_2\h^{\alpha+\kappa}$.\\
Finally,
$$\|w(t)\|_{H^k(M^d)}\lesssim \h^{-(\beta+1)k}\|w(t)\|_{H^k_h}\longrightarrow
0,  $$
for $0\leq t\leq s_2\h^{\alpha+\kappa}$, when $\h\longrightarrow 0$,
what we wanted to prove.

\end{proof}



\section{The instability argument}\label{SectionInstab}
We have now the tools to show our main results.\\
We consider Cauchy conditions $v^0=a^0\e^{iS^0/h}$ of \eqref{nlsh}
which do not oscillate, i.e. such that $S^0=0$. We have seen in the
previous section, that for some analytic amplitudes $a^0$, the
solution writes  $v=a\e^{iS/h}$ and therefore oscillates immediately
with magnitude $\sim\frac1h$.\\[5pt]
Let $\chi$ be given by \eqref{defchi} and $a^0\in
\mathcal{H}^0(l,B)$ nontrivial (for instance $a^0(y)=\e^{-y^2}$).\\
Now set 
\begin{equation}\label{defu0}
u_0^h(x)=\h^{\gamma}\chi(\h^{-\eta}|\kappa(x)|)a^0(\h^{-1}\kappa(x)),
\end{equation}
as initial data for \eqref{nls}.\\
Then we have the Ansatz \eqref{defvapp}, \eqref{defuapp}, \eqref{Muapp}
\begin{equation}\label{Ansatzproof}
u_{\text{app}}(t,x)=\h^{\gamma}\chi(\h^{-\eta}|\kappa(x)|)a(\h^{-\alpha}t,\h^{-1}\kappa(x))\e^{iS(\h^{-\alpha}t,\h^{-1}\kappa(x))/h},
\end{equation}
with $u_{\text{app}}(0,\cdot)=u_0^h$.\\
For $0<c_0\ll 1$ satisfying Proposition \ref{corapproxsc}, set
$$h=\h^{\beta}=\frac{c_0}{n}$$
 with $n\in\N$, and this induces the sequences in
the statements of our main results. In particular 
\begin{equation*}
\text{supp\;}u^h_{0}\subset \{(t,x)\in\R\times M^d\;:\;\;|\kappa(x)|\leq r\h^{\eta}\},
\end{equation*} 
and hence we can choose
$$r_n=\max_{|x|\leq r\h^{\eta}}{|\kappa^{-1}(x)|_g}\longrightarrow 0, $$
in Theorems \ref{thm1} and \ref{thm2}. Here we have assumed $m=0$,
reduction which is always possible.
\subsection{Proof of Theorem \ref{thm1}}
$~$\\[5pt]
\noindent Let $0<\eps<1$ and define 
\begin{equation}\label{defdelta}
\delta_h=\h^{\eps\beta}\log{\frac1{\h}}=\frac1{\beta}h^{\eps}\log{\frac1{h}},
\end{equation}
which tends to 0 with $h$. This choice will become clear later. Consider 
\begin{equation}\label{defu0tilde}
\widetilde{u_0}^h=(1+\delta_h)u_0^h,
\end{equation}
and the associate function $\widetilde{u_{\text{app}}}$.\\[5pt]
In all this subsection we take
$$\gamma =-\frac{d}{p+1}.$$
 This is the right parameter $\gamma$ so that $u_{\text{app}}$ and $\widetilde{u_{\text{app}}}$
are normalized in $L^{p+1}(M^d)$ uniformly for $h\in ]0,\eps[$.

\begin{lemm}\label{lemCI} Let $p\geq (d+2)/(d-2)$ be an odd integer, and let $u_0^h,\widetilde{u_0}^h$ be
  defined by \eqref{defu0}, \eqref{defu0tilde}. Then 
\begin{equation*}
H^+(u_0^h)\lesssim1, \quad H^+(\widetilde{u_0}^h)\lesssim1.
\end{equation*}
There exist $\nu_0>0$ and $q_0>p+1$, such that for all
$0<\nu<\nu_0$ and $p+1\leq q<q_0$
\begin{equation}\label{st5}
\|u_0^h-\widetilde{u_0}^h\|_{H^{1+\nu}(M^d)},\quad
\|u_0^h-\widetilde{u_0}^h\|_{L^{q}(M^d)}\longrightarrow 0,
\end{equation}
when $h\longrightarrow 0$.

\end{lemm}

\begin{proof}
We make the change of variables  $y=\h^{-1}\kappa(x)$, then 
\begin{eqnarray*}
\|\nabla u_0^h\|_{L^2(M^d)}^2&\sim&
\h^{2\gamma+d-2}\int_{|y|\leq
  r\h^{-1+\eta}}\big|\nabla\big(\chi(\h^{1-\eta}y)a^0(y)\big)\big|^2\text{d}y\\
&\sim&
\h^{2\gamma+d-2}\int \big|\nabla a^0(y)\big|^2\text{d}y,
\end{eqnarray*}
as $0<\eta<1$.\\
As $2\gamma +d-2=-2d/(p+1)+d-2>0$ when $p>2d/(d-2)-1$, it follows that 
$$\|\nabla u_0^h\|_{L^2(M^d)}\longrightarrow 0\quad \text{for} \quad h\longrightarrow 0.$$
Compute
$$\| u_0^h\|_{L^{p+1}(M^d)}^{p+1}\sim
\h^{(p+1)\gamma+d}\int |a^0(y)|^{p+1}\text{d}y\sim
\h^{(p+1)\gamma+d}.$$
By definition ${(p+1)\gamma+d}=0$, hence $\|
u_0^h\|_{L^{p+1}(M^d)}$ remains bounded when $h$ tends to $0$, as
well as $H^+(u_0^h)$.\\
Similarly, $ H^+(\widetilde{u_0}^h)\lesssim1$.\\
By the definition \eqref{defdelta} of $\delta_h$ we also have for all $\sigma \geq 0$
\begin{equation}\label{eq345}
\|u_0^h-\widetilde{u_0}^h\|^2_{H^{\sigma}(M^d)}\sim \h^{2\gamma+d-2\sigma}
\delta^2_h\sim \h^{2\gamma+d-2\sigma+2\eps\beta}\big(\log{\frac1{\h}}\big)^2. 
\end{equation}
The terms in \eqref{eq345} tend to 0 if 
$$\sigma<\gamma+\frac{d}2+\eps\beta.$$
But, by \eqref{parameters} and as $p>(d+2)/(d-2)$, 
$$\gamma+\frac{d}2>1,$$
hence we can choose $\nu_0=\eps\beta$ in the statement.\\   
The proof of the other part is similar.
\end{proof}

\begin{proof}[Proof of Theorem \ref{thm1}]
The statements \eqref{st1}, \eqref{st2} and \eqref{st4}
have already been proved in Lemma \ref{lemCI}. \\
Let $0<\eps<1$ which appears in \eqref{defdelta}, and set $s_h=h^{1-\eps}={\h}^{\beta(1-\eps)}$ and $t_h=\h^{\alpha}s_h=\h^{\alpha+\beta(1-\eps) }$.
Denote by $S=S(\h^{-\alpha}t_h,\h^{-1}\kappa(x))$ and by
$b=\chi(\h^{-\eta}|\kappa(x)|)a(\h^{-\alpha}t_h,\h^{-1}\kappa(x))$. Then we have
\begin{eqnarray}
\|(u_{\text{app}}-\widetilde{u_{\text{app}}})(t_h)\|_{L^{p+1}(M^d)}=
\h^{\gamma}  \| b\, \e^{iS/h} -\widetilde{b}\, \e^{i\widetilde{S}/h}\|_{L^{p+1}(\mathcal{U})}\nonumber\\
\geq\h^{\gamma}  \| b\,( \e^{i(\widetilde{S}-S)/h } -1)\|_{L^{p+1}(\mathcal{U})}-
\h^{\gamma}  \| b -\widetilde{b} \|_{L^{p+1}(\mathcal{U})}.\label{eq36}
\end{eqnarray}
We now estimate the l.h.s. terms of \eqref{eq36}.\\
First compute 
\begin{equation*}
\h^{\gamma}  \| (b -\widetilde{b})(\h^{-\alpha}t_h,\h^{-1}\kappa(\cdot)) \|_{L^{p+1}(\mathcal{U})}
\sim \h^{\gamma+d/(p+1)}  \| (b -\widetilde{b})(s_h,\cdot) \|_{L^{p+1}(\mathcal{V})}.
\end{equation*}
From the well-posedness of \eqref{systgeom}, we deduce

\begin{equation*}
 \| (b -\widetilde{b})(s_h,\cdot)
 \|_{L^{p+1}(\h^{-1}\mathcal{V})}\longrightarrow 0,
\end{equation*}
where $s_2$ is given by Proposition \ref{propapprox}.\\
Hence 
\begin{equation}\label{eq40}
\h^{\gamma}  \| (b
-\widetilde{b})(\h^{-\alpha}t_h,\h^{-1}\kappa(\cdot))
\|_{L^{p+1}(\mathcal{U})}\longrightarrow0,\quad h\longrightarrow0.
\end{equation}
Secondly, a  Taylor expansion near $s=0$ shows that 
\begin{equation}\label{taylor}
(\widetilde{S}-S)(s_h,y)\sim -\omega (p-1)\delta_h s_h\big(\chi(\h^{-1-\eta}|y|)a^0(y)\big)^{p-1}.
\end{equation}
Now observe that 
$$ \frac{\delta_hs_h}{h}\sim\log{\frac1h}\longrightarrow +\infty.$$
We then deduce from \eqref{taylor} that for all $|y|\leq 1$,
\begin{equation*}
 \limsup_{h\to 0}\big|\e^{i(\widetilde{S}-S)/h )}-1\big|=2,
\end{equation*}
and as  $\h^{\gamma} \| b(s_h,\h^{-1}\kappa(x))\|_{L^{p+1}(\mathcal{U})} \sim 1$, we obtain 
\begin{equation}\label{eq41}
\limsup_{h\to 0} \h^{\gamma}  \| b\,( \e^{i(\widetilde{S}-S)/h }
-1)\|_{L^{p+1}(\mathcal{U})}\geq c.
\end{equation}
Thus, according to  \eqref{eq40} and \eqref{eq41}
\begin{equation}\label{eq45}
 \limsup_{h\to
 0}\|(u_{\text{app}}-\widetilde{u_{\text{app}}})(t_h)\|_{L^{p+1}(M^d)}\geq c.
\end{equation}

\noindent Finally, if we denote by $L^{p+1}=L^{p+1}(M^d)$,
\begin{eqnarray}
\|(u-\widetilde{u})(t_h)\|_{L^{p+1}}\geq
\|(u_{\text{app}}-\widetilde{u_{\text{app}}})(t_h)\|_{L^{p+1}}-\|(u-{u_{\text{app}}})(t_h)\|_{L^{p+1}}\nonumber\\
-\|(\widetilde{u}
-\widetilde{u_{\text{app}}})(t_h)\|_{L^{p+1}}.\label{eq42}
\end{eqnarray}
If $\eps>0$ is chosen small enough, we can apply  Corollary \ref{Mcorapprox}, with $\kappa=(1-\eps)\beta$, with yields
$\|(u-{u_{\text{app}}})(t_h)\|_{L^{p+1}},\|(\widetilde{u}-\widetilde{u_{\text{app}}})(t_h)\|_{L^{p+1}}
\longrightarrow 0$ with $h$, and thus from \eqref{eq41} and \eqref{eq42}
$$\limsup_{h\to 0}\|(u-\widetilde{u})(t_h)\|_{L^{p+1}}\geq  \limsup_{h\to 0}\|(u_{\text{app}}-\widetilde{u_{\text{app}}})(t_h)\|_{L^{p+1}}>c, $$
which concludes the proof.
\end{proof}
\subsection{Proof of Theorem \ref{thm2}}
$~$\\[5pt]
Here we deal with the case $(M^d,g)=(\R^d,\text{can})$.\\
\noindent Let $\beta>0$ and $\gamma(p-1)=-2(\beta+1)$ as prescribed by
\eqref{parameters}. Let
$0<\sigma<d/2-2/(p-1)$ and 
$$\frac{\sigma}{\frac{p-1}{2}(\frac{d}2-\sigma)}<\rho\leq\sigma.$$
Consider $u_{\text{app}}$ defined by \eqref{Ansatzproof} and let
$s_2>0$ be given by Proposition \ref{propapprox}.\\
Then, according to Corollary \ref{corapprox}, the solution $u$ of
\eqref{nls} with initial condition $u(0)=u_{\text{app}}(0)$ satisfies
for all $k\in \R$
$$
\|(u-u_{\text{app}})(t_h)\|_{H^k(\R^d)}\longrightarrow0,\quad
h\longrightarrow0,  $$
with $t_h=\h^{\alpha }s_2$.
To prove that $u$ satisfies \eqref{1.6} and \eqref{1.7} we only have
to check that 
\begin{equation*}
\|u_{\text{app}}(0)\|_{H^{\sigma}(\R^d)}\to 0 \quad\text{and}\quad
\|u_{\text{app}}(t_h)\|_{H^{\rho}(\R^d)}\to +\infty.
\end{equation*}
To begin with, 
\begin{equation}\label{X}
\|u_{\text{app}}(0)\|_{H^{\sigma}(\R^d)}\sim \h^{\gamma
  -\sigma+d/2}.
\end{equation}
Then, use the equations \eqref{systgeom} to observe that $a\nabla S
 (s_2,\cdot)\not\equiv0$. Hence
\begin{equation}\label{XX}
\|u_{\text{app}}(t_h)\|_{H^{\sigma}(\R^d)}\sim \h^{\gamma-(\beta+1)\rho+d/2}.
\end{equation}
By \eqref{X} and \eqref{XX}, we only have to show that we can choose
$\beta>0$ so that 

\begin{eqnarray}
\gamma -\sigma+d/2&>&0,\label{50}\\
\gamma -(\beta+1)\rho+d/2&<&0\label{51}.
\end{eqnarray}
Let $\eps>0$ such that 
\begin{eqnarray}
\sigma  &<&d/2-2/(p-1)-\eps,\label{52}\\
\rho&>&\frac{\sigma+\eps}{\frac{p-1}{2}(\frac{d}2-\sigma-\eps)},\label{53}
\end{eqnarray}
and take 
\begin{equation}\label{choice1}
\gamma=\sigma-d/2+\eps.
\end{equation}
Therefore by \eqref{52} and \eqref{choice1} we obtain 
\begin{equation}\label{choice2}
\beta=-\frac{p-1}2\gamma-1=-\frac{p-1}2(\sigma-d/2+\eps)-1>0.
\end{equation}
Moreover, with the choice \eqref{choice1}, inequality \eqref{50} is satisfied.\\
Finally,  using the relations \eqref{choice1} and  \eqref{choice2}, we
deduce that \eqref{51} is equivalent to 
\begin{equation*}
\rho>\frac1{\beta+1}(\gamma+d/2)=\frac{\sigma+\eps}{\frac{p-1}{2}(\frac{d}2-\sigma-\eps)},
\end{equation*}
which is satisfied by \eqref{53}.

\subsection{Proof of Theorem \ref{thm3}}
$~$\\[5pt]
\noindent Assume here that $(M^d,g)$ is an analytic riemannian manifold with an
analytic metric $g$.\\
Consider the function $u_{\text{app}}$ defined by \eqref{Ansatzproof} and let
$s_2>0$ be given by Proposition \ref{propapprox}.\\
Let $\kappa\geq 0$ such that $\beta+\eta-\kappa<1$.
 Denote by $t_h=\h^{\alpha+\kappa }s_2$, then by  Corollary
 \ref{Mcorapprox},
 the solution $u$ of
\eqref{nls} with initial condition $u(0)=u_{\text{app}}(0)$ satisfies
for all $k\in \R$
\begin{equation}\label{CV}
\|(u-u_{\text{app}})(t_h)\|_{H^k(M^d)}\longrightarrow0,\quad
h\longrightarrow0.  
\end{equation}
\\[5pt]
$\bullet$ Let 
\begin{equation}\label{choix1}
\frac{d}2-\frac{4}{p-1}<\sigma<\frac{d}2-\frac{2}{p-1}.
\end{equation}
Choose $\eps>0$ so that 
\begin{equation}\label{choix2}
\sigma <\frac{d}2-\frac{2}{p-1}-\eps, 
\end{equation}
and define 
\begin{equation*}
\gamma=\sigma-d/2+\eps,
\end{equation*}
thus
\begin{equation*}
\beta=-\frac{p-1}2\gamma-1=-\frac{p-1}2(\sigma-d/2+\eps)-1.
\end{equation*}
Then by \eqref{choix1} and \eqref{choix2}, $0<\beta<1$. Choose now
$\eta>0$ so small that $0<\beta+\eta<1$. The convergence \eqref{CV}
then follows with $t_h=\h^{\alpha }s_2$.\\
Finally
\begin{equation*}
\|u_{\text{app}}(t_h)\|_{H^{\sigma}(\R^d)}\sim
\h^{\gamma-(\beta+1)\rho+d/2}\longrightarrow +\infty,
\end{equation*}
for 
\begin{equation*}
\rho>\frac1{\beta+1}(\gamma+d/2)=\frac{\sigma+\eps}{\frac{p-1}{2}(\frac{d}2-\sigma-\eps)},
\end{equation*}
which was to prove.\\[5pt]
$\bullet$ Assume here that $0<\sigma<\frac{d}2-\frac4{p-1}$. For
$\beta>0$ and $\eps>0$, take $\kappa=\beta-1+2\eps$ and $\eta=\eps$,
so that $\beta+\eta-\kappa=1-\eps<1$. Then \eqref{CV} holds and 
\begin{equation}\label{equiv}
\|u_{\text{app}}(0)\|_{H^{\sigma}(\R^d)}\sim \h^{\gamma
  -\sigma+d/2}\quad \text{and}\quad \|u_{\text{app}}(t_h)\|_{H^{\sigma}(\R^d)}\sim \h^{\gamma-(\beta+1-\kappa)\rho+d/2}.
\end{equation}
Define $\gamma=\sigma-d/2+\eps$, then 
\begin{equation*}
\beta=-\frac{p-1}2\gamma-1=-\frac{p-1}2(\sigma-d/2+\eps)-1>0,
\end{equation*}
and
\begin{equation*}
\|u_{\text{app}}(0)\|_{H^{\sigma}(\R^d)}\longrightarrow
0\quad\text{when}\quad h\longrightarrow0.
\end{equation*}
The second term in \eqref{equiv} tends to $+\infty$ when 
 \begin{equation*}
\rho>-\frac{\gamma+d/2}{\beta+1-\kappa}=\frac{\sigma+\eps}{2(1-\eps)},
\end{equation*}
which concludes the proof, as $\eps>0$ is arbitrary.


\appendix
\section{Appendix}
Here we reproduce a part of the work of P. G\'erard \cite{PGX}.
\begin{prop}\label{propappendix}
Let  $r>0$ be given by \eqref{defr}. Let  $v_{\text{app}}$ be the function  defined by
\eqref{defvapp}, and let $v$ be the  solution of 
\begin{equation*}
\left\{
\begin{aligned}
&i h\partial_sv +h^2\Delta v   = \omega | v|^{p-1} v,\quad (s,z)\in {\R}^{1+d},\\
&v(0,z)=v_{\text{app}}(0,z).
\end{aligned}
\right.
\end{equation*}
Then there exist $s_2>0$, $\lambda>0$ and $\delta_2>0$ such that 
$v$ can be extended to a continuous function on $[0,s_2]$,
holomorphic-valued on $\big\{(|\h z|<r )\cap (|\text{Im}\;z|<\lambda)\big\}$, and so that for all $k\in\N$
\begin{equation}\label{estpropappendix}
\sup_{|s|<s_2}\;\sup_{|\text{Im\,}z|<\lambda}\;
\|(1-h^2\Delta)^{k/2}(v-v_{\text{app}})(s,\cdot+i\text{Im\,}z)\|_{L^2(\mathcal{B}(0,r/\h))}\leq C_k \e^{-\delta_2/h},
\end{equation}
with $C_k>0$.
\end{prop}

\begin{proof}

Let 
\begin{equation*}
 r_h=-i h\partial_sv_{\text{app}}-h^2\Delta v_{\text{app}}
+\omega ( v_{\text{app}} \overline{v_{\text{app}}})^{\frac{p-1}{2}}
v_{\text{app}}.
\end{equation*}
Then $f_1=v-v_{\text{app}}$ satisfies
\begin{equation}\label{eqA1}
i h\partial_sf_1+h^2\Delta f_1=F(s,z,f_1,h)+r_h,
\end{equation}
where $F$ stands for 
\begin{equation*}
F(s,z,f_2,h)=\omega\big(( v_{\text{app}}
+f_2)^{\frac{p+1}{2}}(\overline{v_{\text{app}}}
+\overline{f_2})^{\frac{p-1}{2}}-
 ( v_{\text{app}} \overline{v_{\text{app}}})^{\frac{p-1}{2}}v_{\text{app}}\big),
\end{equation*}
with $\overline{f}(s,z)=\overline{f(\overline{s},\overline{z})}$.\\
We now show, that for $s_2>0$ and  $\lambda>0$ small enough, there exists
a solution $f_1$ of \eqref{eqA1} such that $f(s,z)$ is continuous on
$[0,s_2]$, holomorphic-valued on $|\text{Im\;}z|<\lambda$, and
exponentially decreasing in $h$: There exists $\delta>0$ so that for
all $k\in \N$
$$\sup_{0\leq s\leq s_2}\;\sup_{|\text{Im\,}z|<\lambda}\;
\|(1-h^2\Delta)^{k/2}f_1(s,\cdot+i\text{Im\,}z)\|_{L^2(\mathcal{B}(0,r/\h))}\leq C_k \e^{-\delta/h}. $$
This will be done thanks to a fixed point argument.\\
For $0< \lambda<l$ and $k>d/2$, set 
\begin{equation*}
\| f \|_h=\sup_{|\text{Im\,}z|<\lambda}\;\|(1-h^2\Delta)^{k/2} f(\cdot+i\text{Im\,}z)\|_{L^2{(\mathcal{B}(0,r/\h))}}.
\end{equation*}
Let $s_1>0$ be given by Proposition \ref{Mpropformel}. Let also
$\delta >0$ and   $s_2\in[0,s_1]$. If $f=f(s,z)$ is continuous on
$[0,s_1]$ and analytic on $|\text{Im\;}z|<\lambda$ we set 
\begin{equation*}
N_h(f,s_2)=\sup_{0\leq s \leq s_2} \e^{\delta(1-\frac{s}{2s_2})/h}\|f(s)\|_h.
\end{equation*}
By the Sobolev embeddings, we have
\begin{equation*}
\sup_{s\leq s_2}\;\sup_{|\text{Im\;}z|<\lambda}\;|f(s,z)|\lesssim N_h(f,s_2)\e^{-\frac{\delta}{4h}}.
\end{equation*}
By Proposition \ref{corapproxsc} we can choose $\delta>0$, $\lambda>0$ and $K>0$ so that  
\begin{equation}\label{AA0}
\| f_1(0)\|_h\leq K\e^{-\frac{\delta}{h}},\quad  N_h(r_h,s_1)\leq
K,\quad \text{and}\quad \sup_{s\leq s_1}\;\sup_{|\text{Im\;}z|<\lambda}\;|v_{\text{app}}(s,z)|\leq K.
\end{equation}
Now use that 
\begin{equation}\label{AA1}
\| f\,g \|_h\lesssim h^{-k}\| f\|_h\| g \|_h,\quad \sup_{s\leq s_2}\; 
\|f(s)\|_h\lesssim  \e^{-\frac{\delta}{2h}} N_h(f,s_2),
\end{equation}
to deduce that 
for all $L>0$, there exists $C_L>0$ and $h_L>0$ so that if $h<h_L$
and $N_h(f,s_2)\leq L$, the following estimates hold  
 \begin{equation}\label{AA2}
\|F(s,\cdot,f(s),h)\|_h\leq C_L\|f(s)\|_h,
\end{equation}
and

\begin{equation}\label{contraction}
\|F(s,\cdot,f_1(s),h)-F(s,\cdot,f_2(s),h)\|_h\leq C_L\|f_1(s)-f_2(s)\|_h.
\end{equation}

\noindent Let $L>0$ to be chosen and $h<h_L$. For $f$ such that $N_h(f,s_2)\leq
L$, let $w$ be the solution of 

\begin{equation}\label{nls2}
\left\{
\begin{array}{l}
ih \partial_t w+h^2\Delta w  = F(s,z,f,h)+r_h,\\
w(0)=f_1(0).
\end{array}
\right.
\end{equation}
The usual $L^2$-estimates for the Schr\"odinger equation and
\eqref{AA1} yield

\begin{equation}\label{AA3}
\|w(s)\|_h\leq\| f_1(0)\|_h+\frac1h\int_0^s\big(C_L\|f(\tau)\|_h+\|r_h(\tau)\|_h \big)\text{d}\tau.
\end{equation}
Now use the exponential decrease of the  norms $\|\,\cdot\,\|_h$ to obtain
\begin{equation*}
\frac1h\int_0^s\|f(\tau)\|_h\text{d}\tau\leq \frac{2s_2}{\eps} N_h(f,s_2) \e^{\delta(1-\frac{s}{2s_2})/h}.
\end{equation*}
Then, with \eqref{AA0}, we deduce from \eqref{AA3}
\begin{equation}\label{AA4}
N_h(w,s_2)\leq K(1+\frac{2s_2}{\eps})+\frac{2s_2}{\eps}C_LN_h(f,s_2).
\end{equation}
Now fix $L\geq 4K$ and $s_2\leq \min{(s_1,\eps/2,\eps/(4C_L))}$, by
\eqref{AA4} we obtain $N_h(w,s_2)\leq L$.\\
We have proved that the application $f\mapsto w$, induced by the
equation \eqref{nls2}, maps the ball
$\{f,\;N_h(f,s_2)\leq L\}$ into itself. Finally use
\eqref{contraction} to show that this application is a
contraction. Hence \eqref{eqA1} admits a unique solution $f_1$,
which satisfies the estimate \eqref{estpropappendix}.
\end{proof}







\end{document}